\documentclass{article}

\usepackage[latin1]{inputenc}
\usepackage[T1]{fontenc}
\usepackage{amsmath}
\usepackage{amssymb}
\usepackage{theorem}
\usepackage[french]{babel}

\newtheorem{theo}{Th{\'e}or{\`e}me}
\newtheorem{defi}{D{\'e}finition}
\newtheorem{lem}{Lemme}
\newtheorem{prop}{Proposition}
\newtheorem{exem}{Exemple}
\newtheorem{cor}{Corollaire}
\newtheorem{rem}{Remarque}

\newcommand{\flh[1]}{#1^{\uparrow}}
\newcommand{\flb[1]}{#1^{\downarrow}}
\newcommand{\oct}{\mathbb{O}}
\newcommand{\R}{\mathbb{R}}
\newcommand{\C}{\mathbb{C}}
\newcommand{\h}{\mathbb{H}}
\newcommand{\g}{\mathfrak{g}}
\newcommand{\G}{\mathcal{G}}
\newcommand{\gt}{\tilde{\mathfrak{g}}}
\newcommand{\B}{\mathcal{B}}
\newcommand{\im}{\text{\normalfont Im\,}}
\newcommand\ssi{si, et seulement si,\ }
\author{Idrisse Khemar}
\title{Surfaces isotropes de $\mathbb{O}$ et syst{\`e}mes int{\'e}grables.}
\date{}

\begin{document}
\maketitle

\section*{Introduction}
Dans cet article, nous {\'e}tudions certaines surfaces isotropes de
$\oct=\mathbb{R}^8$. L'id{\'e}e de s'int{\'e}resser {\`a} de telles surfaces
vient de la volont{\'e} de chercher des analogues aux surfaces
lagrangiennes hamiltoniennes stationnaires de $\mathbb{R}^4$, dans
$\mathbb{R}^8$. Ces surfaces de $\R^4$ forment un syst{\`e}me
compl{\`e}tement int{\'e}grable pr{\'e}sentant une structure in{\'e}dite (cf.\cite{HR1})
 et il est naturel d'en chercher des g{\'e}n{\'e}ralisations dans $\R^8$. (cf. \cite{tern}.)\\
Consid{\'e}rons une surface $\Sigma$ lagrangienne de $\R^4$. On peut
localement trouver une param{\'e}trisation conforme de $\Sigma$ par
des coordonn{\'e}es $(u,v)\in\Omega$ , $\Omega$ {\'e}tant un ouvert de
$\R^2$, i.e. une immersion $X:\Omega \to \R^4$ telle que
$$dX=e^f(e_1 du + e_2 dv)$$ avec $(e_1 ,e_2)$ base hermitienne de
$\mathbb{C}^2$ pour tout $(u,v)\in\Omega$ , $ f\in C^\infty (\R^2)$.
L'identification entre $\R^4$ et $\mathbb{C}^2$ est donn{\'e}e par
$(x_1,x_2,x_3,x_4)\mapsto (x_1 +ix_2,x_3+ix_4)$. A la surface
$\Sigma$ est associ{\'e} \emph{l'angle lagrangien} $\beta$ d{\'e}fini par
 $ e^{i\beta}=\det(e_1,e_2)$ (qui ne d{\'e}pend pas de la
 param{\'e}trisation choisie car il ne d{\'e}pend que du plan tangent
 $T_{X(u,v)}\Sigma =\R e_1 +\R e_2 $).\\
 Maintenant consid{\'e}rons la fonctionnelle d'aire
 $\mathcal{A}(\Sigma)=\int_{\Sigma}dv $ sur l'ensemble des
 surfaces orient{\'e}es lagrangiennes de $\R^4$. Un point critique
 pour cette fonctionnelle est une surface lagrangienne $\Sigma$
 telle que $\delta\mathcal{A}(\Sigma)(X)=0$ pour tout champ de
 vecteur $X$ {\`a} support compact sur $\R^4$ avec la condition suppl{\'e}mentaire
 que $X$ doit {\^e}tre lagrangien i.e. son flot preserve les surfaces
 lagrangiennes; si on suppose que cela n'est vrai que si $X$ est hamiltonien,
 i.e. $X=-i\nabla h $ avec $h\in C^\infty(\R^4,\R)$ alors $\Sigma$
 est dite hamiltonienne stationnaire. On montre alors, cf.\cite{HR1},
 que $\Sigma$ est hamiltonienne stationnaire \ssi $\triangle\beta
 =0$ o{\`u} $\triangle$ est le laplacien sur $\Sigma$ d{\'e}fini {\`a} l'aide
 de la m{\'e}trique induite.\\
 Dans \cite{HR1}, il est montr{\'e} que les surfaces lagrangiennes
 hamiltoniennes stationnaires de $\R^4$ sont solutions d'un syst{\`e}me
 compl{\`e}tement int{\'e}grable. Dans le langage des syst{\`e}mes compl{\`e}tement
 int{\'e}grables on construit une famille de connexions de courbure nulle
 $\alpha_{\lambda}$, $\lambda\in S^1$ qui s'{\'e}crit:
 \begin{equation}\label{alphalambda}
  \alpha_{\lambda}=\lambda^{-2}\alpha'_2 +\lambda^{-1}\alpha'_{-1}
 +\alpha_0 +\lambda\alpha''_1 +\lambda^2 \alpha''_2
 \end{equation}
 cf. \cite{HR1}.\\
 Nous nous proposons ici de trouver des surfaces de $\R^8$ isotropes
 telles qu'{\`a} chacunes d'elle, $\Sigma$, corresponde une fonction
 $\rho_\Sigma :\Sigma \to S^3 $ (analogue de
 $e^{i\beta}:\Sigma\to S^1$) et que les surfaces pour lesquelles
 $\rho_\Sigma$ est harmonique forment un syst{\`e}me compl{\`e}tement
 int{\'e}grable. Pour ce faire nous allons proc{\'e}der par analogie avec les
 surfaces lagrangiennes hamiltoniennes stationnaires de $\R^4$.
 Commen{\c c}ons par formuler le probl{\`e}me dans $\R^4$ {\`a} l'aide des
 quaternions. Ensuite nous proc{\`e}derons par analogie en utilisant les
 octonions (section~2). Le rappel des d{\'e}finitions et connaissances n{\'e}cessaires
 sur les octonions est fait dans la section~1.\\
 Pour $x,y\in \mathbb{H}$, on a
 $$ x\cdot\bar{y}=\langle x,y\rangle_{\R^4} -\,\omega (x,y)i\,
 -{\det}_{\mathbb{C}^2}(x,y)j=\langle x,y\rangle_{\mathbb{C}^2} -
 {\det}_{\C^2}\,(x,y)j$$ o{\`u} $\omega=dx_1\wedge dx_2 + dx_3\wedge dx_4
 $. Ainsi pour $(e_1,e_2)$ base orthonorm{\'e}e d'un plan lagrangien de
 $\R^4$ on a
 $$e_1\cdot\bar{e}_2 =-e^{i\beta}j$$
o{\`u} $\beta$ est l'angle lagrangien du plan $\text{Vect}(e_1 ,e_2)$. On
 voit que l'on peut exprimer la contrainte lagrangienne ainsi
que l'angle lagrangien {\`a} l'aide du produit dans $\h$. Plus
pr{\'e}cisemment, nous avons appliqu{\'e} la proc{\'e}dure suivante: on forme le produit
 $x\cdot \bar{y}$, avec $x,y$ de
norme 1, puis on suppose que les deux premiers termes  de
$x\cdot\bar{y}\in\h$ dans la d{\'e}composition $\h=\R\oplus\R
i\oplus\R j\oplus\R k$, sont nuls et alors on a
$x\cdot\bar{y}\in S^1 j $ ce qui nous permet de recup{\'e}rer
$e^{i\beta}\in S^1 $.\\
Nous proc{\`e}derons de m{\^e}me dans $\R^8 $ en utilisant le produit des
octonions (section~2). Pour $q,q'\in\oct=\h^2$ on a
$$q\cdot\overline{q'}=\langle q,q'\rangle +\sum_{i=1}^7 \omega_i(q,q')e_i$$ o{\`u}
$(e_i)_{0\leq i\leq 7}$ est la base canonique de $\R^8$,
 $\omega_i=\langle \cdot,L_{e_i}\cdot\rangle $, et $L_{e_i}$ designe la
  multiplication {\`a} gauche par $e_i$. Nous introduisons la d{\'e}composition
$$q.\overline{q'}=(B(q,q'),-\rho(q,q'))\in\h^2.$$
 Alors nous regardons les
surfaces isotropes pour $\omega_1 ,\,\omega_2,\,\omega_3$, que nous
appelons \emph{surfaces} $\Sigma_V$. Nous nous interessons donc {\`a}
l'ensemble $Q$ des plans de $\oct$ isotropes pour ces 3 formes
symplectiques. L'ensemble $V$ des bases orthonorm{\'e}es de ces plans est
l'ensemble des couples norm{\'e}es $(q,q')\in S^7\times S^7 $ qui
v{\'e}rifient $B(q,q')=0$ (V est l'analogue de l'ensemble des bases hermitiennes de
$\C^2$). On a alors $Q=V/SO(2)$, et $\rho(q,q')\in S^3$ pour tout
$(q,q')\in V$ (la norme  est multiplicative dans $\oct$ :
$|qq'|=|q||q'|=1)$. On a ainsi d{\'e}fini une fonction
$\rho:V\to S^3$ analogue {\`a} $e^{i\beta}$ et cette fonction passe au
quotient en une application $\rho:Q\to S^3$.\\ D'autre part dans le  cas de
$\R^4$ on a le groupe $U(2)$ qui agit librement et transitivement
sur l'ensemble des bases hermitiennes et on peut {\'e}crire l'action
de $U(2)$ {\`a} l'aide des quaternions: on a le morphisme surjectif de
groupe de noyau $\pm 1$:
 $$\begin{array}{rcl}
  S^3\times S^3 & \to & SO(4)\\
  (p,q) & \mapsto & L_p R_{\bar{q}}=R_{\bar{q}}L_p =(x\mapsto
  px\bar{q}) \end{array}$$
et on peut {\'e}crire
 $$x.\bar{y}=\langle x,y\rangle +\langle x,iy\rangle i+\langle x,jy\rangle j
  + \langle x,ky\rangle k$$
Ainsi $U(2)$ est le sous-groupe de $SO(4)$ qui commute avec $L_i$
d'o{\`u} $U(2)=\{L_p R_{\bar{q}}\, ,\, p\in S^1 ,\, q\in S^3 \}$. Quant
{\`a} $SU(2$) c'est le sous-groupe de $SO(4)$ qui commute avec
$L_i,L_j,L_k $ d'o{\`u} $SU(2)=\{R_q\, , \, q\in S^3 \}$.
D'une mani{\`e}re g{\'e}n{\'e}rale, pour $g=L_p R_{\bar{q}}\in SO(4)$ on a
$(gx)(\overline{gy})=p(x\bar{y})\bar{p}$.\\
Par analogie nous allons chercher le
groupe qui conserve $B$, i.e. le sous-groupe de $SO(8)$ qui
conserve $\omega_1,\,\omega_2,\,\omega_3$. Nous trouverons le groupe
$SU(2)\times SU(2)$. Le r{\'e}sultat est tr{\`e}s diff{\'e}rent de ce qui ce
passe dans $\h$ mais c'est le bon groupe de sym{\'e}trie. En effet
comme $\rho_{\Sigma}$ est {\`a} valeurs dans $S^3$ et que nous voulons lui
imposer d'{\^e}tre harmonique, nous allons utiliser la th{\'e}orie des
applications harmoniques du point de vue des syst{\`e}mes int{\'e}grables
(cf. \cite{DPW}) et donc {\'e}crire $S^3$ comme un espace
sym{\'e}trique: $S^3= S^3\times S^3 / \triangle$ o{\`u} $\triangle$ est la
diagonale de $S^3\times S^3$, $\triangle=\{(a,a)\, ,\, a\in S^3\}$
(et est l'analogue de $SU(2)$). Cependant nous nourrissons l'espoir
d'avoir un groupe qui agit transitivement sur $V$ (ce qui n'est
pas le cas de $S^3\times S^3$) ainsi qu'il en est pour $U(2)$ et
les bases hermitiennes de $\C^2$. Nous allons donc chercher {\`a} grossir le
groupe en prenant le sous-groupe de $SO(8)$ qui conserve la
nullit{\'e} de $B$: nous trouvons alors un groupe $G$ de dimension 9 ($V$
est de dimension 10) qui n'agit donc pas transitivement. Alors nous
 regardons l'action de $G$ sur $Q$, qui lui est de
dimension 9. Malheureusement nous trouvons que l'action n'est
toujours pas transitive. Nous calculons alors les orbites: nous
trouvons que toutes les orbites sont de dimension 8 sauf deux orbites
d{\'e}g{\'e}n{\'e}r{\'e}es l'une de dimension 7, l'autre de dimension 6. En outre
nous construisons une fonction $p\colon Q\to [ 0\, ,\, 1/2 ]$ dont les
fibres sont les orbites de $G$, les orbites d{\'e}g{\'e}n{\'e}r{\'e}es {\'e}tant $p^{-1}(\{0\})$
et $p^{-1}(\{ \frac{1}{2}\})$ respectivement. Ensuite nous arrivons {\`a}
trouver un moyen simple de passer d'une orbite {\`a} une autre (cf. th{\'e}or{\`e}me
\ref{alphabeta}). Enfin nous terminons la section~2 en {\'e}tudiant
alg{\'e}briquement le groupe $G$ et son alg{\`e}bre de Lie.\\
Dans la section~3 nous {\'e}tudions les surfaces $\Sigma_V $. Nous montrons que
celle dont le $\rho_{\Sigma}$ est harmonique forment un syst{\`e}me
compl{\`e}tement int{\'e}grable: nous contruisons une famille de connexion de
courbure nulle $\alpha_{\lambda}$ comme dans (\ref{alphalambda}).\\
Ensuite nous montrons que les surfaces $\Sigma_V$ sont solutions de
deux {\'e}quations l'une lin{\'e}aire, l'autre non lin{\'e}aire. (En fait c'est
la m{\^e}me {\'e}quation o{\`u} on repr{\'e}sente la surface de mani{\`e}re
diff{\'e}rente).\\
Dans la section~4, nous exposons la m{\'e}thode des groupes de lacets, en
se r{\'e}f{\'e}rant {\`a} \cite{DPW} et \cite{HR1} pour les d{\'e}tails.
Puis nous obtenons une
repr{\'e}sentation de type Weierstrass pour les surfaces $\Sigma_V$.\\
Dans la section~5, nous calculons le vecteur courbure moyenne d'une
surface $\Sigma_V$ (dans l'espoir d'obtenir une interpr{\'e}tation
variationnelle).\\
Dans la section 6, nous montrons que ce que nous avons fait pour les surfaces
$\Sigma_V$ est en fait un cas particulier de quelque chose de plus
g{\'e}n{\'e}ral. En effet, en consid{\'e}rant le produit vectoriel de $\oct$, nous
d{\'e}finissons  une application $\rho\colon Gr_2(\oct)\to S^6$. Alors
nous montrons que les surfaces  immerg{\'e}es $\Sigma$ de $\oct$ telles que
$\rho_{\Sigma}\colon z\in \Sigma \mapsto\rho(T_z\Sigma)\in S^6$ est
harmonique (surfaces $\rho-harmoniques$) forment un syst{\`e}me
compl{\`e}tement int{\'e}grable. Le groupe de sym{\'e}trie est alors $Spin(7)$.
Plus g{\'e}n{\'e}ralement, soit $I\varsubsetneqq\{1,...,7\}$ alors les
surfaces $\omega_I-isotropes$, i.e. isotropes pour $\omega_i$,
$i\in I$, dont le $\rho_{\Sigma}$ (qui est alors {\`a} valeurs dans
$S^{I}=S(\oplus_{i \notin I,i>0}\R e_i)\simeq S^{6-|I|}$) est
harmonique, forment un syst{\`e}me compl{\`e}tement int{\'e}grable. Le groupe de
sym{\'e}trie est alors $G_I\simeq Spin(7-|I|)$. Pour $I=\{1,2,3\}$, on
retrouve les surfaces $\Sigma_V$. Nous  construisons donc une famille
$(\mathcal{S}_I)$ param{\'e}tr{\'e}e par $I$, d'ensembles de surfaces solutions d'un
syst{\`e}me int{\'e}grable, tous inclus
dans $\mathcal{S}_\emptyset$, telle que $I\subset J$ implique $\mathcal{S}_J
\subset\mathcal{S}_I$. Par restriction {\`a} $\h$ on obtient les surfaces
$\rho-harmoniques$, $\omega_I-isotropes$ de $\h$. Alors $\rho(Gr_2(\h))=S^2$
 et $|I|=0,1$ ou $2$. Pour $|I|=1$ on retrouve les surfaces lagrangiennes
hamiltoniennes stationnaires de $\R^4$ et pour $|I|=2$, les surfaces
sp{\'e}ciales lagrangiennes. Par restriction {\`a} $\text{Im}(\h)$, on
retrouve les surfaces CMC de $\R^3$. Nous terminons l'article par le
calcul du vecteur courbure moyenne d'une surface quelconque de
$\oct$ en fonction de $\rho$ (nous pensons qu'il existe une interpr{\'e}tation
variationnelle des surfaces $\rho -\text{harmoniques}$).
%
%
%
%
%

  \section{L'alg{\`e}bre des octonions}
\subsection{D{\'e}finitions}

On rappelle ici les d{\'e}finitions et propri{\'e}t{\'e}s sur les octonions
qui nous seront utiles pour la suite. Pour avoir plus de d{\'e}tails et
pour les d{\'e}monstrations on pourra consulter \cite{Har},\cite{HaL}.
On appelle alg{\`e}bre des octonions l'espace vectoriel:
$$
 \oct = \left\{ \left(\begin{array}{cc}
     x & -\bar{y} \\
     y & \bar{x}
   \end{array} \right) \, ,\, x,y\in\h \right\}\subset M_2(\h)
$$
muni de la multiplication (i.e. application bilin{\'e}aire sur $\oct$):
$$
\left( \begin{array}{cc}
  x & -\bar{y} \\
  y & \bar{x}
\end{array} \right)\cdot \left( \begin{array}{cc}
  x' & -\bar{y}' \\
  y' & \bar{x}'
\end{array} \right) = \left( \begin{array}{cc}
  xx'-y'\bar{y} & -\bar{y}\bar{x}'-\bar{y}'x \\
  x'y +\bar{x}y' & \bar{x}'\bar{x}-y\bar{y}'
  \end{array} \right)
$$
On voit que l'on peut identifier $\oct$ {\`a} $\h^2=\R^8$ muni de la
multiplication
$$(x,y)\cdot(x',y')=(xx'-y'\bar{y},\, x'y+\bar{x}y')$$
On d{\'e}finit une conjugaison $q\mapsto\bar{q}$ {\`a} l'aide de la
conjugaison sur les matrices: $\overline{(x,y)}=(\bar{x},-y)$. On
remarque la pr{\'e}sence d'un {\'e}l{\'e}ment neutre pour la multiplication:
$\mathbf{1}=(1,0)$ et $\oct$ est ainsi une $\R$-alg{\`e}bre
unitaire. En particulier, $\R .\mathbf{1}$ est une sous-alg{\`e}bre
dont les {\'e}l{\'e}ments seront dits r{\'e}els et caract{\'e}ris{\'e}s par
 $\bar{q}=q$.\\
On d{\'e}finit sur $\oct$ la norme
$N(q)=q.\bar{q}=\bar{q}.q=x\bar{x}+ y\bar{y}\in \R$ qui n'est
autre que la norme euclidienne standard de $\h^2=\R^8$. \emph{Cette
norme est multiplicative}: $N(qq')=N(q)N(q')$. En particulier $S^7$ est
stable par multiplication. Les octonions orthogonaux {\`a} $\mathbf{1}$
, $\R^{\perp}$ , pour la norme $N$ seront dits : \emph{octonions
pures}, et caract{\'e}ris{\'e}s par $\bar{q} = -q$ , ou encore $q^2 \in
\R^-$. Si on se restreint {\`a} la sph{\`e}re $S^7 =\{q\in \oct \, ,\,
\bar{q}.q =1\}$ alors le sous-ensemble des octonions purs de
$S^7$ est caract{\'e}ris{\'e} par $q^2=-1$.\\
La base canonique de $\R^8$ correspond en {\'e}criture matricielle {\`a}
la base des octonions:
$$
 \flh[E]=
\left(\begin{array}{cc}  1 & 0 \\  0 & 1 \end{array}\right),\:
 \flh[I]=
\left(\begin{array}{cc}  i & 0 \\  0 & -i \end{array}\right),\:
 \flh[J]=
\left(\begin{array}{cc}  j & 0 \\  0 & -j \end{array}\right),\:
 \flh[K]=
\left(\begin{array}{cc}  k & 0 \\  0 & -k \end{array}\right)
$$
$$
 \flb[E]=
\left(\begin{array}{cc}  0 & -1 \\  1 & 0 \end{array}\right),\:
 \flb[I]=\flb[E]\flh[I]=
\left(\begin{array}{cc}  0 & i \\  i & 0 \end{array}\right),\:
 \flb[J]=\flb[E]\flh[J]=
\left(\begin{array}{cc}  0 & j \\  j & 0 \end{array}\right),
$$
$$
\flb[K]=\flb[E]\flh[K]=
\left(\begin{array}{cc}  0 & k \\  k & 0 \end{array}\right).
$$
Dans la suite il nous arrivera aussi de noter cette base $(e_i)_{0\leq
i\leq 7}$ (l'ordre des vecteurs {\'e}tant toujours le m{\^e}me).\\
$\oct$ n'est pas associative: $\flb[I](\flb[J]\flb[K])=\flb[I](-\flh[I])=\flb[E] $
tandis que $(\flb[I]\flb[J])\flb[K]=(-\flh[K])\flb[K]=-\flb[E]$.

\subsection{Propri{\'e}t{\'e}s de la multiplication}

\begin{prop}\
\begin{description}
\item[(i)]  $\langle xz,yz\rangle =N(z)\langle x,y\rangle,\
\langle zx,zy\rangle =N(z)\langle x,y\rangle$
\item[(ii)]  $\langle xz,yw\rangle + \langle yz,xw\rangle =
2\langle x,y\rangle \langle z,w\rangle$
\item[(iii)]  si $x\notin \R 1 $, $\R 1\oplus\R x$ est une alg{\`e}bre
isomorphe (isom{\'e}trique) {\`a} $\C$.
\end{description}
\end{prop}

\begin{prop}\
\begin{description}
\item[(i)] $\overline{xy}=\bar{y}\bar{x}$ , $ \bar{\bar{x}}=x$ ,
$\langle x,y\rangle =Re(x\bar{y})=Re(\bar{x}y)$
\item[(ii)] $x(\bar{x}y)=N(x)y$ , $(x\bar{y})y=N(y)x$
\item[(iii)] $x(\bar{y}z) + y(\bar{x}z)=2\langle x,y\rangle z$ , $(z\bar{y})x
+ (z\bar{x})y =2\langle x,y\rangle z$
\item[(iv)] si $\langle x,y\rangle =0$ , alors $x\bar{y}=-y\bar{x}$ et
$x(\bar{y}z)=-y(\bar{x}z)$ , $(z\bar{y})x=-(z\bar{x})y$.
\end{description}
\end{prop}

\begin{defi} On dira d'un {\'e}l{\'e}ment $x\neq 0$ d'une alg{\`e}bre $A$ non
associative qu'il est inversible s'il existe $x'\in A$ tel que $\forall y \in
A$ , $x'(xy)=x(x'y)=(yx)x'=(yx')x=y$. Ceci revient {\`a} dire que $
L_x \colon y\mapsto xy$ et $R_x \colon y\mapsto yx$ sont
inversibles d'inverses respectives $L_{x'}$ et $R_{x'}$.
\end{defi}

\begin{prop} Tout $x\in \oct\smallsetminus\{0\}$ est inversible
d'inverse $N(x)^{-1}\bar{x} $. En outre on a $^t L_x=L_{\bar{x}}$ ,
$^t R_x=R_{\bar{x}}$ .
\end{prop}
\begin{prop}On a les propri{\'e}t{\'e}s d'associativit{\'e} suivantes:
\begin{description}
\item[(i)] \begin{eqnarray*}(ax)(ya) & = & a((xy)a) \\
                           a(x(ay)) & = & (a(xa))y \\
                           x(a(ya)) & = & ((xa)y)a
            \end{eqnarray*}
\item[(ii)] \begin{eqnarray*} (xy)x & = & x(yx)\\
                              x(xy) & = & x^2y \\
                              (xy)y & = & xy^2
            \end{eqnarray*}
\item[(iii)] $R_xL_x=L_xR_x$ , $L_x^2 =L_{x^2}$ , $R_x^2=R_{x^2}$,
$L_{axa} =L_a L_x L_a $ , $R_{aya}=R_aR_yR_a$.

\end{description}
\end{prop}

\begin{prop} L'application trilin{\'e}aire $\{x,y,z\}=(xy)z-x(yz)$ est
altern{\'e}e donc antisym{\'e}trique.
\end{prop}

\begin{prop}\
\begin{description}
\item[(i)] $x,y\in \oct$ commutent \ssi $(1,x,y)$ est li{\'e}e et alors
la sous-alg{\`e}bre (unitaire) engendr{\'e}e par $x$ et $y$ est isomorphe {\`a}
$\C$ (on suppose  que $\{x,y\}\nsubseteq \R 1$).
\item[(ii)] s'ils ne commutent pas alors la sous-alg{\`e}bre qu'il
engendrent est isomorphe (isom{\'e}trique) au corps $\h$ des
quaternions.
\item[(iii)] si $D$ est une sous-alg{\`e}bre (unitaire) de $\oct$ de
dimension 4 (i.e. $\simeq\h$) et si $a\in D^{\perp}\smallsetminus\{0\}$
alors $\oct=D\oplus a.D$ et on a
$$ (x +ay)(x' + ay')=(xx' - \lambda y'\bar{y}) + a(x'y + \bar{x}y')$$
avec $\lambda = -N(a)$.
\item[(iv)] soit $a,b \in (\R 1)^{\perp}$ unitaires et orthogonaux
(alors comme $ab=-ba$) la sous-alg{\`e}bre $D$ engendr{\'e}e par $a,b$ est
de dimension 4, et soit $c\in D^{\perp}$ unitaire alors
$\{1,a,b,ab,c,ca,cb,c(ab)\}$ est une base orthonorm{\'e}e de $\oct$
qui a la m{\^e}me table de multiplication que la base canonique.
\end{description}
\end{prop}

\begin{prop}
 Si $x,y$ ne commutent pas i.e. $(1,x,y)$ est libre
alors $z\mapsto \{x,y,z\}$ n'est pas identiquement nulle autrement
dit $L_{xy} \neq L_xL_y$ ou encore par antisym{\'e}trie $R_{yx} \neq
R_xR_y$ , ou encore $L_xR_y \neq R_yL_x$.
\end{prop}

\begin{theo}
Si $L_xL_y=L_z$ alors obligatoirement $z=xy$ et donc $x,y$
commutent, et de m{\^e}me $R_xR_y =R_z\Longrightarrow z=yx$. Ainsi
$\{L_x ,\ x\in \oct^*\}$ et $\{ R_x ,\ x\in \oct^*\}$ ne sont pas
des sous-groupes de $GL(\oct)=GL(\R^8)$.
\end{theo}

\begin{prop}
$L_xL_y =L_yL_x\Longleftrightarrow xy=yx $ (idem pour $R$).
\end{prop}

%
%

\section{Plans isotropes et Groupes op{\'e}rants}
\subsection{ Analogie {\`a} l'aide des octonions}

Comme nous l'avons expliqu{\'e} dans l'introduction, on va {\'e}tudier
l'expression $q.\overline{q'}$ dans $\oct$ par analogie avec
$\h$. Soit donc $q=(x,y),\, q'=(x',y') \in \oct$ , alors on a
$$ q.\overline{q'}={x\choose y} .{\overline{x'}\choose {-y'}}=
{{x\overline{x'} + y'.\overline{y}}\choose {\overline{x'}y -
\overline{x}y'}}.$$
On peut aussi {\'e}crire en notant $(e_i)_{0\leq i\leq 7}$ la base
canonique de $\R^8=\oct$ d{\'e}finie {\`a} la section~1:
$$q\cdot \overline{q'}=\langle q,q'\rangle + \sum_{i=1}^7
\langle q,\ e_i\cdot q'\rangle e_i \ .$$
Posons $\omega_i (q,q')=\langle q,\ e_i\cdot q'\rangle$ , $0\leq i\leq 7$,
alors $\omega_i$ est la forme symplectique sur $\R^8$ associ{\'e}e {\`a}
 l'endomorphisme $L_{e_i}$ ($(L_{e_i})^2=-Id$). On a alors
$$ q\cdot \overline{q'}=\langle q,q'\rangle +
\sum_{i=1}^7 \omega_i(q,q')\, e_i \ .$$
On notera
$$B(q,q')=x\,\overline{x'} + y'\overline{y}=\langle q,q'\rangle +
\sum_{i=1}^3 \omega_i (q,q')\, e_i$$
et
$$\rho(q,q')=\overline{x}y' - \overline{x'}y=
-\sum_{i=4}^7 \omega_i(q,q')\, e_{i-4} \ .$$
Soit alors $V=\{(q,q')\in S^7 \times S^7 / \ B(q,q')=0\}$. On a alors
 pour tout $(q,q')\in V$
 $$ q\cdot \overline{q'}={0 \choose {-\rho}}$$
 avec $\rho\in S^3$. Ceci s'{\'e}crit encore
  $$q'= {0\choose \rho}\cdot q .$$
 On a ainsi d{\'e}fini une fonction $\rho\colon (q,q')\in V \mapsto\rho(q,q') \in S^3$.
 On peut calculer les coordonn{\'e}es de $q'$ en fonction de celles de
 $q$ d'apr{\`e}s l'expression pr{\'e}c{\'e}dente :
 $$ q'=\left(\begin{array}{c}
   -y\bar{\rho} \\
    x\rho   \
 \end{array}\right).
 $$
En particulier on voit que $V$ est une sous vari{\'e}t{\'e} de $S^7\times
S^7$ diff{\'e}omorphe {\`a} $S^7\times S^3$, le diff{\'e}omorphisme {\'e}tant
{\'e}videmment $(q,q')\mapsto (q,\rho)$. Enfin on remarque que $\rho$
ne d{\'e}pend que du plan orient{\'e} engendr{\'e} par $(q,q')$.

\subsection{A la recherche de groupes agissant sur $V$}

Cherchons le sous-groupe de $GL(8)$ qui conserve $B$, i.e. le groupe
des {\'e}l{\'e}ments $g\in SO(8)$ qui commutent avec
$L_{\flh[I]},L_{\flh[J]} ,L_{\flh[K]}$ , i.e. avec les
$L_{(x,0)}$, $x\in \h$. On a $L_{(x,0)}=\left(\begin{array}{cc}
  L_x & 0 \\
  0 & L_{\bar{x}}
\end{array}\right)$ ainsi en posant $g=\left(\begin{array}{cc}
  A & B \\
  C & D
\end{array}\right)$ on a pour tout $x\in \h$ :
$$\left( \begin{array}{cc}  L_x & 0 \\  0 & L_{\bar{x}} \end{array}\right)
\left(\begin{array}{cc}  A & B \\  C & D \end{array}\right) -
 \left(\begin{array}{cc}   A & B \\  C & D  \end{array}\right)
\left( \begin{array}{cc}    L_x & 0 \\  0 & L_{\bar{x}}   \end{array}\right)
=$$
$$\left(\begin{array}{cc}
  L_xA - AL_x & L_xB - BL_{\bar{x}} \\
  L_{\bar{x}}C - CL_x & L_{\bar{x}}D - DL_{\bar{x}}
\end{array}\right).$$
En {\'e}galant la derni{\`e}re matrice {\`a} 0 on obtient: $\forall x\in\h \
[L_x,A]=[L_x,D]=0$ , $BL_x=L_{\bar{x}}B ,\ CL_x=L_{\bar{x}}$. Les
{\'e}quations sur $A,D$ signifient que $A=R_a$, $D=R_d$  avec $a,d
\in\h$. Pour $B,C$ on a $B(x.1)=\bar{x}B(1)$ d'o{\`u}
$B(ax)=\overline{ax}B(1)$ mais on a aussi $B(ax)=B(L_ax)=\bar{a}B(x)
=\bar{a}\bar{x}B(1)$ or comme $\overline{ax}\neq\bar{a}\bar{x}$ en
g{\'e}n{\'e}ral on a donc $B(1)=0$ et donc $B=0$ et de m{\^e}me pour $C$. D'o{\`u}

\begin{theo}
$g\in SO(8)$ conserve $B$ \ssi
$$g=\left(\begin{array}{cc}
  R_a & 0 \\
  0 & R_d
\end{array}\right)$$
avec $a,d\in S^3$, ainsi le groupe conservant $B$ est $SU(2)\times SU(2)$.
\end{theo}

\begin{rem}\emph{ Le groupe obtenu n'agit pas transitivement. D'autre
part, contrairement {\`a} ce qui se passe dans $\h$, ici, comme on l'a
vu dans la section~1, $R_q$ et $L_p$ ne commutent pas, $\{R_qL_p
,\ q,p\in S^7\} $ n'est pas un groupe et on n'a pas:
$(qb)\overline{(q'b)}=q\overline{q'}$.\\ On peut se demander
pourquoi le r{\'e}sultat est tr{\`e}s diff{\'e}rent de ce qui se passe avec
les quaternions o{\`u} le groupe qui conserve $\langle\cdot ,\cdot\rangle_{\C}$, i.e. $U(2)$, agit
transitivement sur les bases hermitiennes. Cela s'explique par l'absence d'associativit{\'e}.
Dans $\h$, si $g\in SO(4)$ commute avec $L_i$ et $L_j$, il commute
alors avec $L_k=L_iL_j$ tandis que dans $\oct$ le fait de commuter
avec $L_{\flh[K]}$ est une condition suppl{\'e}mentaire. D'ailleurs, on
peut voir que le sous-groupe de $SO(8)$ qui commute avec
$L_{\flh[I]},L_{\flh[J]}$ est de dimension 10 (car isomorphe {\`a}
$Sp(2)=U(2,\h)$) et donc la condition de commutativit{\'e} avec
$L_{\flh[K]}$, fait passer la dimension de 10 {\`a} 6 . Le groupe,
$S^3\times S^3$, ainsi obtenu est le bon groupe de sym{\'e}trie
recherch{\'e} pour obtenir un syst{\`e}me compl{\`e}tement int{\'e}grable, car
pour {\'e}crire $S^3$ sous la forme d'un espace sym{\'e}trique , on peut
{\'e}crire $S^3=S^3\times S^3/\triangle$ o{\`u} $\triangle$ est la
diagonale.\\ De m{\^e}me dans \cite{HR1}, il suffirait de se restreindre au
groupe $S^1$, or  les auteurs y utilisent le groupe $U(2)$
tout entier qui agit transitivement sur les bases hermitiennes, en
{\'e}crivant que $S^1=U(2)/SU(2)$. On voudrait ici de la  m{\^e}me
mani{\`e}re trouver un groupe (contenant $S^3\times S^3$) qui
agisse transitivement (ou du moins "le plus transitivement
possible") sur $V$. Car dans le cas o{\`u} $G$ agit transitivement sur
$V$ on peut toujours relever un couple $(q,q')\in V$ en un {\'e}l{\'e}ment
de $G$, et le fait de travailler sur $G$ permet de repr{\'e}senter une
surface par un {\'e}l{\'e}ment de $\R_{+}^{*}.G$ (dans \cite{HR2} on obtient
ainsi une {\'e}quation  de  Dirac). D'autre part, on veut comprendre
la g{\'e}om{\'e}trie de $V$, i.e., la g{\'e}om{\'e}trie des plans isotropes pour les
trois formes symplectiques  $\omega_i$, $i=1,2,3$ de la m{\^e}me fa{\c c}on que la
g{\'e}om{\'e}trie des plans lagrangiens de $\C^2$ est {\'e}tudi{\'e}e (ou rappel{\'e}e) dans
\cite{HR1}.}
\end{rem}

Le plus gros groupe qui conserve $V$ est le sous groupe de $SO(8)$
qui conserve la nullit{\'e} de $B$. Il est donn{\'e} par :

\begin{theo}
Soit $G=\{g\in SO(8)/ B(q,q')=0\Longleftrightarrow B(g.q,g.q')=0\}$,
c'est le plus grand sous-groupe de $GL(8)$ qui stabilise $V$. Il
existe un morphisme surjectif $\theta\colon G \to O(3)\subset O(4)$
{\`a} valeurs dans le groupe des isom{\'e}tries de $\h$ qui fixent 1, tel que
$\theta^{-1}(SO(3))=G^0$, la composante neutre de $G$, et
$\theta^{-1}(O^{-}(3))=L_{\flb[E]}G^0$, donc
$$G=G^0\sqcup L_{\flb[E]}G^0 \ ,$$
 et tel que tout $g\in G^0$
s'{\'e}crit $g=\left(\begin{array}{cc}
  R_a\theta(g) & 0 \\
  0 & R_b\theta(g)
\end{array}\right)$.
Plus Pr{\'e}cisement
$$G^0=\left\{\left(\begin{array}{cc}
  R_aL_c & 0 \\
  0 & R_bL_c
\end{array}\right),\ a,b,c\in S^3 \right\}$$
De plus on a $\forall q,q'\in\oct$, $B(g.q,g.q')=\theta(g)(B(q,q'))$
pour tout $g\in G$. En outre pour tout $g\in G^0$, on a
 $\theta(L_{\flb[E]}g)=\theta(g)\ast=\ast\, \theta(g) $
 o{\`u} $\ast$ d{\'e}signe la conjuguaison de $\h$ qui
est aussi l'{\'e}l{\'e}ment $-I_3$ de $O(\im\h)=O(3)$. Enfin dans $G^0$, on a
$\theta(Diag(R_aL_c,R_bL_c))=Int_c=(x\in \im\h \mapsto cxc^{-1})$.
\end{theo}

\noindent {\em D{\'e}monstration} --- Dire que $g\in GL(8)$ v{\'e}rifie
$B(q,q')=0\Longrightarrow B(g.q,g.q')=0$ revient {\`a} dire que
$(\langle q,L_{e_i}q'\rangle =0,\ 0\leq i\leq 3)\Longrightarrow
(\langle g.q,L_{e_i}g.q'\rangle =0, \ 0\leq i\leq 3)$ ce qui {\'e}quivaut {\`a}
\begin{equation}\label{theta-ij}
^t g L_{e_i}g=\sum_{j=0}^3 \theta_{ij}L_{e_i},\quad 0\leq
i\leq 3.
\end{equation}
avec $\theta_{ij}\in \R$. Alors en posant
$\theta(g)=(\theta_{ij})_{0\leq i,j\leq 3}$, on a
$\theta(gg')=\theta(g)\theta(g')$ pour tout $g,g'\in G'=\{g\in
GL(8)/\ B(q,q')=0\Longleftrightarrow B(g.q,g.q')=0\}$. En effet on a
\begin{eqnarray*}
^t(gg') L_{e_i}(gg') = {}^tg'(\sum_{j=0}^3
\theta_{ij}(g)L_{e_j})g' & = & \sum_{j=0}^3\sum_{k=0}^3
\theta_{ij}(g)\theta_{jk}(g')L_{e_k} \\
 & = & \sum_{k=0}^3(\theta(g)\theta(g'))_{ik}L_{e_k}
\end{eqnarray*}
d'o{\`u} le r{\'e}sultat. Ainsi $\theta\colon G'\mapsto GL(4)$ est un
morphisme de groupe. Prenons $i=0$ dans (\ref{theta-ij}),
comparons l'{\'e}quation obtenue avec sa transpos{\'e}e, alors en utilisant
que $ ^t L_{e_j}=-L_{e_j}$ pour $j\geq 1$, on voit que l'on
doit avoir $^t gg =\theta_{00}Id$ et $\theta_{0j}=0 $
pour $1\leq j\leq 3$. En proc{\'e}dant de m{\^e}me pour $i\geq 1$ on
obtient $\theta_{i\,0}=0$ pour $1\leq i \leq 3$ ; il en r{\'e}sulte que
$\theta=Diag(\theta_{00},\ \mu)$ avec $\mu\in Gl(3)$. De plus comme
$^t gg=\theta_{00}Id$, on doit avoir $\theta_{00}>0$ et
donc $G'=\R_{+}^*.G$ avec, rappelons le, $G=G'\bigcap SO(8)$.
Maintenant en {\'e}crivant (\ref{theta-ij}) pour $i\geq 1$ et $g\in
G$, et en utilisant le fait que les $L_{e_i},\ i\geq 1$,
anticommutent deux {\`a} deux, on obtient :
$$
\begin{array}{ccl}
 -2\delta_{ik}Id & = & g^{-1}(L_{e_i}L_{e_k}+L_{e_k}L_{e_i})g \\
 & = & (g^{-1}L_{e_i}g)(g^{-1}L_{e_k}g) + (g^{-1}L_{e_k}g)(g^{-1}L_{e_i}g) \\
     &  = &  \displaystyle\sum_{1\leq j,l \leq 3}\mu_{ij}\mu_{kl}
    (L_{e_j}L_{e_l} + L_{e_l}L_{e_j})\\
     &  = & \displaystyle\sum_{j=1}^3 \mu_{ij}\mu_{kj}(-2Id)
\end{array}
$$
d'o{\`u} $\mu(g)\in O(3)$ si $g\in G$.\\
Cherchons, ensuite {\`a} quelles conditions $g=\begin{pmatrix}
  A & B \\
  C & D
\end{pmatrix}\in G$. Pour cela , on utilise toujours
(\ref{theta-ij}), et le fait que $L_{e_i}=\left(\begin{array}{cc}
  L_{e'_i} & 0 \\
  0 & -L_{e'_i}
\end{array}\right)$ avec $(e'_1,e'_2,e'_3)=(i,j,k)$
 la base canonique de $\im\h$, cela donne:
$$\begin{pmatrix}
 L_{e'_i}A & L_{e'_i}B \\
  -L_{e'_i}C & -L_{e'_i}D
\end{pmatrix}=
\begin{pmatrix}
 AL_{\mu^{i}} & -BL_{\mu^{i}} \\
 CL_{\mu^{i}}  & -DL_{\mu^{i}}
\end{pmatrix}$$
o{\`u} $\mu^{i}=\sum_{j=0}^{3}\mu_{ij}e'_j\in \R i\oplus\R
j\oplus\R k=\im\h$ .
Ainsi pour $A$ par exemple, on a $L_{e'_i}A=AL_{\mu^{i}}$, d'o{\`u}
$e'_i.A(1)=A(\mu^{i})$ ce qui imlique que
$$A.\left(\begin{array}{cc}
  1 & \begin{array}{ccc}
    0 & 0 & 0
  \end{array} \\
  \begin{array}{c}
    0 \\
    0 \\
    0
  \end{array} & ^t \mu
\end{array}\right)
=(a,\ e'_1.a,\ e'_2.a,\ e'_3.a)=R_a$$
avec $a=A(1)$. Finalement on a donc $A=R_a.Diag(1,\ \mu)=R_a.\theta$
 et de m{\^e}me $D=R_d.Diag(1,\ \mu),\ B=R_b.Diag(1,-\mu),\
C=R_c.Diag(1,-\mu)=R_c.Diag(1,\ \mu)\ast$ , o{\`u} $\ast=Diag(1,-I_3)$
 est la conjuguaison dans $\h$.\\
Ensuite on {\'e}crit qu'on doit avoir $L_{e'_i}.A(e'_j) =
A(\mu^{i}.e'_j)$ pour $1\leq i,j\leq 3$ en utilisant l'expression
de $A$ que l'on vient d'obtenir. On trouve, apr{\`e}s calcul, que:
$$(L_{e'_i}=AL_{\mu^{i}},\ 1\leq i\leq 3\,) \Longleftrightarrow
(\mu=com(\mu)\ \text{ou}\ a=0\,)$$
($com$ d{\'e}signe la comatrice), comme $\mu\in O(3)$ cela veut dire
$\det\mu =1$ ou $a=0$. On trouve la m{\^e}me chose pour $D$.
Pour $B$ on a aussi:
$$(L_{e'_i}=-BL_{\mu^{i}},\ 1\leq i\leq 3\,) \Longleftrightarrow
(\det(\mu)=-1\ \text{ou}\ b=0\,)$$
 et de m{\^e}me pour $C$. On ach{\`e}ve la d{\'e}monstration en remarquant que
 $L_{\flb[E]}=\begin{pmatrix}
   0 & -Id \\
   Id & 0 \
 \end{pmatrix}$, et que le groupe des automorphismes de $\h$ est {\'e}gal au
 groupe des automorphismes int{\'e}rieurs de $\h$ qui n'est autre que
$SO(\im\h)$.\hfill $\blacksquare$ \\

\noindent Ce th{\'e}or{\`e}me permet de voir comment $\rho\colon
(q,q')\mapsto \rho(q,q')=\bar{x}y'-\bar{x'}y$ se transforme sous
l'action de $G$: $$\rho(g.q,g.q')=\bar{a}\rho(q,q')b \quad
\text{pour} \quad g=Diag(R_aL_c,\ R_bL_c)\ .$$ Ainsi l'action de
$G^0$ sur $V\cong \{(q,\rho)\in S^7\times S^3\}$ s'{\'e}crit $g\cdot
(q,\rho)= (g.q,\,\bar{a}\rho\, b)$. En outre, on voit que l'action de
$G$ sur $\rho$ d{\'e}finit une action transitive de $G^0$ sur $S^3$.
Si l'on oublie les $L_c$ qui n'ont aucun effet sur $\rho$, et que
l'on se restreint au groupe $S^3\times S^3$, cette action n'est
autre que le rev{\^e}tement universel de $SO(4)$. Maintenant pour
$g'=Diag(R_aL_c,R_bL_c).L_{\flb[E]}$ on a
$$\rho(g'.q,g'.q')=\bar{a}\,\overline{\rho(q,q')}\,b  \ .$$ On a
donc trouv{\'e} un groupe $G$ de dimension 9 agissant sur $V$ qui est
de dimension 10. Cette action ne peut donc pas {\^e}tre transitive. On
est donc amen{\'e} {\`a} {\'e}tudier l'action de $G$ sur les plans engendr{\'e}s
par les {\'e}l{\'e}ments de $V$, en esp{\'e}rant qu'elle soit transitive.

\subsection{Action de $G$ sur les plans de $V/SO(2)$}

Consid{\'e}rons l'ensemble des plans orient{\'e}s de $\R^8$ engendr{\'e}s par les
 $(q,q')\in V$:
 $$Q=\{\text{Plans orient{\'e}s annulant }\omega_1,\,\omega_2,\,\omega_3\}.$$
$Q$ est une sous-vari{\'e}t{\'e} compact de $Gr_2(\R^8)=\{\text{plans
orient{\'e}s de }\R^8\}$, diff{\'e}omorphe {\`a} $V/SO(2)$, l'action de $SO(2)$ sur
$V$ {\'e}tant donn{\'e}e par $(q,q')\mapsto (q,q')\cdot R_{\theta}$\hspace{0cm}$=(
\cos\theta\,q + \sin\theta\,q',-\sin\theta\,q + \cos\theta\,q')$.
En effet, comme $V$ est ferm{\'e} dans $Y=\{(q,q')\text{ orthonorm{\'e}es
de }\R^8\}$, le graphe $R_V$ de l'action de $SO(2)$ sur $V$, $R_V=
R_Y \cap(V\times V)$, est ferm{\'e}, puisque $R_Y$, le graphe de
$Y$ modulo $SO(2)$ est ferm{\'e} ($Y/SO(2)=Gr_2(\R^8)$ a une structure
de vari{\'e}t{\'e} quotient) donc $V/SO(2)$ admet une structure de vari{\'e}t{\'e}
quotient et munie de cette structure c'est une sous-vari{\'e}t{\'e} de
$Gr_2(\R^8)$ de dimension 9.\\
Lorsqu'on identifie $V$ {\`a} $S^7\times S^3$ l'action de $SO(2)$
s'{\'e}crit
$$
 R_{\theta}\cdot(q,\rho)=(\cos\theta\,q + \sin\theta  \begin{pmatrix}
  0 \\
  \rho
\end{pmatrix}\cdot q,\ \rho\,) \ .
$$
$G$ agit sur $Q$: si on note $[q,q']$ le plan engendr{\'e} par $(q,q')$ on a
$g.[q,q']=[g.q,g.q']$. Cette action n'est pas transitive. En effet,
 consid{\'e}rons $P_0=[1,\flb[E]]$ alors
 $ G^0.P_0=\{[\genfrac(){0pt}{1pt}{ca}{0},\genfrac(){0pt}{1pt}{0}{cb}]
 ,a,b,c \in S^3 \}=\{[\genfrac(){0pt}{1pt}{x}{0},\genfrac(){0pt}{1pt}{0}{y}
 ],|x|=|y| =1\}$ est de dimension au plus 6 < 9. De
 plus $L_{\flb[E]}[1,\flb[E]]=[\flb[E],-1]=[1,\flb[E]]$. Finalement
$G.P_0=G^0.P_0$ est de dimension 6.\\

\subsubsection{Calcul du stabilisateur d'un point}

Soit $P_0=[q_0,q'_0]\in Q$ et $Stab(P_0)$ le stabilisateur de $P_0$
, alors comme $G$ est compact , on sait que $O(P_0)=G.P_0$ est une
sous-vari{\'e}t{\'e} compacte de $Q$ et $G/Stab(P_0)\cong G.P_0$. Ainsi si
on avait $\dim(Stab(P_0))=0$, alors $O(P_0)$ serait une sous vari{\'e}t{\'e}
de dimension {\'e}gale {\`a} $\dim(G)=\dim(Q)$, donc un ouvert de $Q$ mais
aussi un ferm{\'e} car elle est compact, donc comme $Q$ est connexe (car
$V \simeq S^7 \times S^3$ est connexe) on aurait que $O(P_0)=Q$ donc
$G$ agirait transitivement or ce n'est pas le cas. Donc $\forall P_0
\in Q,\,\dim(Stab(P_0))>0$ et $\dim\,O(P_0)\leq 8$.\\

\noindent Les orbites sont en fait donn{\'e}es par :

\begin{theo}
$Q$ admet la partition suivante:
$$
Q=G.P_1\sqcup G.P_2\sqcup U
$$
o{\`u} \begin{description}
\item[$\bullet$] $P_1=[1,\flb[E]]$ , $ P_2=\left[\begin{pmatrix} \frac{1}{\sqrt{2}}
  \\ \frac{i}{\sqrt{2}}  \end{pmatrix},
 \begin{pmatrix} \frac{-i}{\sqrt{2}}  \\  \frac{1}{\sqrt{2}}  \end{pmatrix}
 \right]$
\item[$\bullet$] $U=\{P=\left[\begin{pmatrix}  x \\  y \ \end{pmatrix},
\begin{pmatrix}  x' \\  y' \ \end{pmatrix}\right]\in Q,/ (x,x') \text{ est
 libre et } |\langle x,x'\rangle | + ||x| - |x'|| \neq 0 \}$
 est un ouvert de $Q$.
\end{description}
De plus on a $G.P_1=\{P\in Q /(x,x') \text{ est li{\'e}e}\}$ , $G.P_2=\{P\in Q
/\langle x,x'\rangle =0 ,\, |x|=|x'| \}$ et enfin $\forall P\in U$, $G.P$ est
une sous-vari{\'e}t{\'e} compacte de dimension 8. Ainsi, il y a deux orbites
d{\'e}g{\'e}n{\'e}r{\'e}es $G.P_1$, et $G.P_2$, et toutes les autres orbites sont de
dimensions 8.\\
On peut ajouter que $\forall P\in Q$, $G^0.P=G.P$ et que
\begin{itemize}
\item[$\bullet$] si $P\in G.P_1$, alors $Stab(P)_{\mid P}=\{\pm Id_P \}$
\item[$\bullet$] si $P\in GP_2$, alors $Stab(P)_{\mid P}=SO(P)$
\item[$\bullet$] si $P\in U$, alors $Stab(P)_{\mid P}=\{\pm Id_P\}$.
\end{itemize}
\end{theo}

\noindent {\em D{\'e}monstration} --- Dans un premier temps, on se
restreint {\`a} l'action de $G^0$. Soit donc $P=[q,q']\in Q$ et posons
$q=(x,y)$, $q'=(x',y')=(-y\bar{\rho},\,x\rho)$. Soit $g\in Stab(P)$,
alors il existe $\theta\in \R$ tel que
$$
g.(q,q')=(q,q')\cdot R_{\theta} \quad i.e. \quad g_{\mid P}=R_{\theta}
$$
ce qui s'{\'e}crit
 \begin{equation}\label{x,x'}
 \begin{array}{lcr}
 \left\{ \begin{array}{l}
          cx\,a = \cos\theta \,x + \sin\theta \,x'\\
          cx'a = -\sin\theta \,x + \cos\theta \,x'
          \end{array} \right.
 & et & \bar{a}\rho \,b = \rho.
 \end{array}
 \end{equation}
Si $(x,x')$ est libre alors cela implique que
$\text{Mat}_{(x,x')}({R_aL_c}_{\mid [x,x']})=R_{\theta}$, ce qui
n{\'e}cessite que ou bien $\langle x,x'\rangle =|x'|-|x|=0$ ou bien $\theta=0$
 mod $\pi$.
 Ceci nous am{\`e}ne {\`a} diff{\'e}rencier trois cas:
\begin{enumerate}
\item $(x,x')$ est libre et ($\langle x,x'\rangle\neq 0$  ou  $|x'| - |x| \neq 0$ )
\item $(x,x')$ est libre et $\langle x,x'\rangle =|x'| - |x| =0$
\item $(x,x')$ est li{\'e}e.
\end{enumerate}
Dans chaque cas, on d{\'e}termine le sous-groupe de $SO(4)$,
$\{R_aL_c\in SO(4)$ v{\'e}rifiant les 2 premi{\`e}res {\'e}quations de (\ref{x,x'})$\}$,
 ce qui nous donne alors $\dim Stab(P)$ et $Stab(P)_{\mid P}$. On trouve
alors que $\dim Stab(P)=1,2,3$ dans les cas 1,2,3 respectivement.
Ceci nous donne $\dim O(P)$ dans chaque cas. Ensuite, dans les cas 2 et 3,
respectivement, on d{\'e}termine facilement un {\'e}l{\'e}ment $g=Diag(R_aL_c,R_bL_c)\in
G^0$ tel que $P=g.P_2$ et $P=g.P_1$ respectivement. Enfin, on
v{\'e}rifie que $L_{\flb[E]}P_1=P_1$, $L_{\flb[E]}P_2=P_2$ et que $\forall
P\in U,\, L_{\flb[E]}P\in G^0.P$ d'o{\`u} $G.P=G^0.P$, $\forall P\in
Q$. Ceci ach{\`e}ve la d{\'e}monstration.\hfill $\blacksquare$ \\

\subsubsection{Caract{\'e}risation des orbites}

On cherche une fonction $p\colon Q\to\R$ dont les fibres soient les
orbites de $G^0$ (et donc de $G$). Elle est donn{\'e}e par le :\\

\begin{theo}\label{p}
Soit $p\colon [q,q']\in Q \mapsto |\im\!(x.\bar{x'})| \in
[0,\frac{1}{2}]$. Alors les orbites de $G$ sont les fibres de $p$:
\begin{enumerate}
\item $p^{-1}({0})=G.P_1=\{P\in Q /(x,x') \text{ est li{\'e}e}\}$
\item $p^{-1}({\frac{1}{2}})=G.P_2=\{P\in Q /\langle x,x'\rangle = |x'| - |x|=0\}$
\item $p^{-1}(]0,\frac{1}{2}[)=U$ et $\forall P,P'\in U$, $p(P)=p(P')
\Longleftrightarrow G.P=G.P'$.
\end{enumerate}
\end{theo}

\noindent {\em D{\'e}monstration} --- D'abord $p$ est bien d{\'e}finie car
$\im\!(x.\bar{x'})$ ne d{\'e}pend que du plan $[q,q']$. Ensuite, elle est
bien {\`a} valeurs dans $[0,\frac{1}{2}]$ puisque $ |x'|^2 + |x|^2 =1$.
 Elle est invariante sous l'action de $G^0$ : pour
$g=Diag(R_aL_c,R_bL_c)$ on a $p([g.q,g.q'])=|c\im\!(x.\bar{x'})c^{-1}|=
|\im\!(x.\bar{x'})|=p([q,q'])$, et sous l'action de $L_{\flb[E]}$ on a
$p(L_{\flb[E]}[q,q'])=\im\!(-y.\overline{(-y')})|=|\im\!(x.\bar{x'})|$.\\
Montrons r{\'e}ciproquement que toute fibre est incluse dans une orbite. On a
$p([q,q'])=0 \Longleftrightarrow x.\bar{x'}=\alpha\in \R\Longleftrightarrow
|x'|^2x=\alpha x'\Longleftrightarrow (x,x')$ est li{\'e}e, d'o{\`u}
$p^{-1}({0})=G.P_1$. Pour la suite on aura besoin du lemme suivant:

\begin{lem}
Pour tout $P\in Q$, il existe un repr{\'e}sentant $(q,q')\in Q$ tel que
\hbox{$\langle x,x'\rangle =0$}.
\end{lem}
{\em D{\'e}monstration du lemme} --- On suppose que $\langle x,x'\rangle \neq 0$
alors on a
\begin{eqnarray*}
\langle\cos\theta\,x + \sin\theta\,x',-\sin\theta\,x + \cos\theta\,x'\rangle & =
& \cos(2\theta)\langle x,x'\rangle \\
 & & + \frac{|x'|^2 - |x|^2}{2}\sin(2\theta)\\
& = & A\cos(2\theta + \phi)
\end{eqnarray*}
et cette derni{\`e}re fonction de $\theta$ s'annule pour certaines
valeurs de $\theta$, ce qui veut dire qu'il existe un repr{\'e}sentant
de $P$ tel que $\langle x,x'\rangle =0$.\\

\noindent {\em D{\'e}monstration du th{\'e}or{\`e}me} --- Dor{\'e}navant, on suppose que
$\langle x,x'\rangle =0$ \ et donc $\im\!(x.\bar{x'})=x.\bar{x'}.$ Ainsi si
$|x.\bar{x'}|=\frac{1}{2}$, alors $|x||x'|=\frac{1}{2}$ et comme
$|x|^2 + |x'|^2=1 $ on a donc $|x|=|x'|=\frac{1}{\sqrt{2}}$, ainsi
comme $\langle x,x'\rangle =0$, on a $P\in G.P_2$. On a bien
$p^{-1}(\{\frac{1}{2}\})=G.P_2$.\\
Soit $P,P'\in U$ tel que $p(P)=p(P')$ avec $P=[q,q'],\,P'=[q_1,q'_1]$.
Alors $p(P)=p(P')\Longleftrightarrow \exists\mu\in SO(3)$ tel que
$\im\!(x_1.\bar{x_1'})=\mu(\im\!(x.\bar{x'}))=c\im\!(x.\bar{x'})c^{-1}$, alors
quitte {\`a} remplacer $(q,q')$ par $(g.q,g.q')$, avec $g=Diag(L_c,L_c)$,
 on est ramen{\'e} au cas o{\`u}
$$
\im\!(x_1.\bar{x_1'})=\im\!(x.\bar{x'})
$$
et en prenant des repr{\'e}sentants convenables, ceci s'{\'e}crit encore
$$
x_1'.\bar{x_1'}=x.\bar{x'}
$$ i.e.
\begin{equation}\label{rho}
x_1\rho_1\,\bar{y_1}=x\rho\,\bar{y}
\end{equation}
d'o{\`u}
$$
\begin{array}{ccccc}
  \left\{\begin{array}{l}
 |x||y|=|x_1||y_1| \\
 |x|^2 + |y|^2=|x_1|^2 + |y_1|^2
\end{array}\right. & \Longrightarrow & \left\{\begin{array}{l}
 |x|=|x_1| \\
 |y|=|y_1|
\end{array}\right. & ou & \left\{\begin{array}{l}
 |x|=|y_1| \\
 |y|=|x_1|
 \end{array}\right.
\end{array}.
$$
On peut se ramener {\`a} la premi{\`e}re des 2 possibilit{\'e}es quitte {\`a} remplacer
$(q,q')$ par $(-q',q)$, on peut alors poser $\left\{\begin{array}{l}
x_1=x.a \\ y_1=y.b \end{array}\right.$, avec $a,b\in S^3$, et en revenant {\`a}
(\ref{rho}), on obtient $\bar{a}\rho\,b=\rho_1$, donc $g.(q,q')=(q_1,q_1')$
 avec $g=Diag(R_a,R_b)$ et finalement
$$
G.P=G.P'
.$$
Ceci ach{\`e}ve la d{\'e}monstration du th{\'e}or{\`e}me.\hfill $\blacksquare$ \\

\begin{theo}\label{alphabeta}
Soit $V_1=\R_{+}^*.\pi_{SO(2)}^{-1}(U)=\{(q,q')\in \oct^*\times
\oct^* /|q|=|q'|,\,B(q,q')=0 \ et \ 0<|\im\!(x.\bar{x'})|<\frac{1}{2}|q|^2
\}$. Alors consid{\'e}rons l'action de $(\R_{+}^*)^2$ sur $V_1$ d{\'e}finie
par :
$$
 (\alpha,\beta)\cdot (q,q')=\left(\begin{pmatrix}
  \alpha x \\
  \beta y
\end{pmatrix},\begin{pmatrix}
  \beta x' \\
  \alpha y'
\end{pmatrix}\right)
$$
Alors l'action de $(\R_{+}^*)^2$ commute avec celle de $G^0$. Soit
$(q_0,q_0')\in V_1$ tel que $\langle x_0,x_0'\rangle =0$ alors
$$
\forall(q,q')\in V_1 ,\,\exists(g,(\alpha,\beta))\in G^0 \times (\R_{+}^*)^2,
\,\exists \theta \in \R \quad tel\ que$$
$$
(g\cdot (\alpha,\beta)\cdot (q_0,q_0'))\cdot R_{\theta}=(q,q')
$$
Il y a exactement deux possibilit{\'e}es pour $(\alpha,\beta)$, l'une tel que
$\alpha < |q|/(\sqrt{2}|x_0|)$, l'autre tel que $\alpha > |q|/(\sqrt{2}|x_0|)$
. $R_{\theta}$ peut {\^e}tre chang{\'e} en $-R_{\theta}$ (et donc $g$ en $-g$), $g$
peut varier dans $g.Stab([q_0,q_0'])$.
\end{theo}

\noindent{\em D{\'e}monstration} --- Posons
$(\tilde{q}_0,\tilde{q}_0')=\frac{1}{|q_0|}(q_0,q_0')$ et $(\tilde{q},\tilde{q}')=
\frac{1}{|q|}(q,q')$.
 Alors on a $p((\alpha',\beta')\cdot
 (\tilde{q}_0,\tilde{q}_0'))=\alpha'\beta' p(\tilde{q}_0,\tilde{q}_0')$.
On choisit $\alpha',\beta'$ tel que $\alpha'^2|\tilde{x}_0|^2 +
\beta'^2|\tilde{y}_0|^2=1$ et
$\alpha'\beta'p([\tilde{q}_0,\tilde{q}_0'])=p([\tilde{q},\tilde{q}'])$,
on montre que ceci est possible et qu'il y a exactement deux solutions
par l'{\'e}tude de la fonction $\alpha'\mapsto
\frac{\alpha'(1-\alpha'^2|\tilde{x}_0|^2)^{\frac{1}{2}}}
{(1- |\tilde{x}_0|^2)^{\frac{1}{2}}}$ dont la valeur maximale est
$\frac{1}{2|\tilde{x}_0||\tilde{y}_0|}=\frac{1}{2p(\tilde{q}_0,\tilde{q}_0')}$ .\\
Ensuite, on pose
$(\alpha,\beta)=(\alpha'\frac{|q|}{|q_0|},\,\beta'\frac{|q|}{|q_0|})$.
Alors $(\alpha,\beta)\cdot (q_0,q_0')$ est de norme $|q|$ et
$[\frac{1}{|q|}(\alpha,\beta)\cdot (q_0,q_0')]$ est dans l'orbite de
$[\tilde{q},\tilde{q}']$ donc il existe $g\in G^0$ et $\theta \in
\R$ tel que $(g\cdot (\alpha,\beta)\cdot (q_0,q_0'))\cdot
R_{\theta}=(q,q')$. \hfill $\blacksquare$ \\

\begin{rem}{\em
L'action de $(\R_{+}^*)^2$ permet de passer d'une orbite {\`a} l'autre
tandis que $G^0$ agit transitivement sur chaque orbite. Cependant,
l'action de $(\R_{+}^*)^2$ n'est pas compatible avec celle de
$SO(2)$: elle n'envoie pas un plan sur un autre plan.
Comme on le voit sur la d{\'e}monstration, le th{\'e}or{\`e}me est valable pour les
{\'e}l{\'e}ments de $G.P_2$, i.e. on aurait pu prendre $ \R_{+}^*.\pi_{SO(2)}^{-1}
(Q\smallsetminus G.P_1)=\{(q,q')\in \oct^*\times\oct^* /|q|=|q'|,\,B(q,q')=0
\ et \ \im\!(x.\bar{x'})\neq 0 \}$ au lieu de $\R_{+}^*.\pi_{SO(2)}^{-1}(U)$.
Dans le cas o{\`u} $[q,q']\in G.P_2$, $(\alpha,\beta)$ est unique et donn{\'e}
par $(\alpha,\beta)=(|q|/(\sqrt{2}|x_0|),\,|q|/(\sqrt{2}|y_0|$ $))$. En ce
qui concerne le point de r{\'e}f{\'e}rence $(q_0,q_0')$, on ne peut pas le
prendre quelconque; comme on le voit sur la d{\'e}monstration on a
besoin de $|\tilde{x}_0||\tilde{y}_0|=p(\tilde{q}_0,\tilde{q}_0')$,
i.e. $\langle x_0,x_0'\rangle =0$. On peut prendre par exemple
$$
(q_0,q_0')=
  \left(\begin{pmatrix}
    \frac{1}{\sqrt{2}}  \\ \frac{i}{\sqrt{2}}
      \end{pmatrix},
  \begin{pmatrix}
   \frac{-i}{\sqrt{2}}  \\  \frac{1}{\sqrt{2}}  \end{pmatrix}\right)
$$
et alors $\rho(q_0,q_0')=1$. On a alors $\alpha^2 + \beta^2=2|q|^2$ et les
deux possibilit{\'e}es pour $(\alpha,\beta)$ dans le cas de $U$  sont
$\alpha <|q|$ et $\alpha >|q|$ tandis que pour $G.P_2$ on a
$\alpha=\beta=|q|$. On voit aussi que le th{\'e}or{\`e}me est encore valable pour
$(q,q')\in G.P_1$ et on alors $\alpha\beta=0$ et $\alpha^2 + \beta^2=2|q|^2$
, mais {\'e}videmment on ne peut pas prendre $(q_0,q_0')\in G.P_1$.}
\end{rem}
\begin{rem}{\em
L'angle $\theta$ a une d{\'e}finition intrins{\`e}que. En effet, soit $P\in
U$, alors il existe un couple $(q_1,q_1')\in P$ unique {\`a} $\pm 1$ pr{\`e}s
tel que $\langle x_1,x_1'\rangle =0$, alors $\theta$ est l'angle d{\'e}fini
 modulo $\pi$ tel que $(q,q')=(q_1,q_1')\cdot R_{\theta}$ pour $[q,q']= P$.
On a donc d{\'e}fini une fonction
$$
(q,q')\in V_1 \mapsto \theta(q,q')\in \R/\pi\mathbb{Z}
$$
Ainsi dans chaque plan $P\in U$, il existe un axe privil{\'e}gi{\'e}. Cette axe
permet de mesurer des angles de droites dans $P$.}
\end{rem}

\subsection{D{\'e}composition de $G$ et de son alg{\`e}bre de Lie.}

On a d{\'e}fini l'application $\rho$ par analogie avec le d{\'e}terminant
sur les bases hermitiennes de $\C^2$. On voudrait d{\'e}finir l'analogue
du d{\'e}terminant sur le groupe $U(2)$, i.e., une fonction d{\'e}finie sur
$G$ {\`a} valeurs dans $S^3$ qui corresponde d'une certaine mani{\`e}re {\`a}
$\rho$.\\
Soit $(q_0,q_0')\in V$ tel que $\rho(q_0,q_0')=1$. Consid{\'e}rons
$\rho(g.q_0,g.q_0')$ pour $g\in G^0$, on a
$$
\rho(g.q_0,g.q_0')=\bar{a}.1.b=\bar{a}.b \ .
$$
On se demande {\`a} quelle condition a-t-on $\rho(g.q_0,g.q_0')=\rho(g'.q_0,g'.q_0')
$. C'est le cas \ssi $a^{-1}.b=a'^{-1}b'\Longleftrightarrow
a'a^{-1}=b'b^{-1}\Longleftrightarrow \exists d\in S^3 /
  \left\{\begin{array}{l}
   b'=db \\
   a'=da
  \end{array}\right.$ i.e. \ $g^{-1}g'=
  \begin{pmatrix}
    R_dL_c & 0 \\
    0 & R_dL_c
  \end{pmatrix}$.
Donc si $g_0\in G^0$ est tel que $\rho(g_0.q_0,g_0.q_0')=\rho_0$
alors $\rho(g.q_0,g.q_0')=\rho_0$ \ssi $g=g_0.\begin{pmatrix}
    R_dL_c & 0 \\
    0 & R_dL_c
  \end{pmatrix}$.
Pour $g_0$, on peut prendre par exemple $g_0=
  \begin{pmatrix}
    Id &  0\\
    0 & R_{\rho_0}
  \end{pmatrix}$.
On a en fait le th{\'e}or{\`e}me suivant :

\begin{theo}\label{semidirect}
Soit
$$G_0^0=\left\{
  \begin{pmatrix}
    R_aL_c & 0 \\
     0 & R_aL_c
  \end{pmatrix},\,a,c\in S^3\right\},\ G_2^0=\left\{ \begin{pmatrix}
    Id &  0\\
    0 & R_{\rho}
  \end{pmatrix},\,\rho\in S^3\right\}$$
Alors:
\begin{description}
  \item[(i)] $G_2^0$ est un sous-groupe distingu{\'e} dans $G^0$: \ $G_2^0
  \vartriangleleft G^0$.
  \item[(ii)] $G^0$ est le produit semi-direct de $G_2^0$ et $G_0^0$:
  \ $G^0=G_2^0\rtimes G_0^0$.
  \item[(iii)] Ceci permet de d{\'e}finir une application
  $\ \tilde{\rho}\colon G^0 \to S^3$ d{\'e}finie par
  $$
\begin{pmatrix}
    Id &  0\\
    0 & R_{\rho}
  \end{pmatrix}.\begin{pmatrix}
    R_aL_c & 0 \\
     0 & R_aL_c
  \end{pmatrix} \mapsto \rho \ .
$$
En outre si $(q_0,q_0')\in V$ est tel que $\rho(q_0,q_0')=1$ alors
$\tilde{\rho}$ est aussi donn{\'e}e par $\tilde{\rho}\colon g\in G^0
\mapsto\rho(g.q_0,g.q_0')$. De plus $\tilde{\rho}$ est invariante
par multiplication {\`a} droite par $G_0^0$:
$\tilde{\rho}(g.h)=\tilde{\rho}(g) \ \forall g\in G^0,\,\forall h\in
G_0^0$.
\item[(iv)] $\rho(g^{-1}\cdot q,\,g^{-1}\cdot q')=1\Longleftrightarrow
\tilde{\rho}(g)=\rho(q,\,q')$, pour $g\in G^0$, $(q,q')\in V$.
\end{description}
\end{theo}

\noindent\emph{D{\'e}monstration} --- Pour (i),(ii) et (iii), c'est un simple
 calcul. Pour (iv) il suffit d'utiliser  $\rho(g.q,g.q')=\bar{a}\rho(q,q')b$.
\hfill $\blacksquare$\\

\begin{rem}{\em
On voit que les $L_c$ ne jouent aucun r{\^o}le dans le th{\'e}or{\`e}me
pr{\'e}c{\'e}dent. Il r{\'e}sulte en particulier de ce dernier que l'on a:
$$
S^3=G^0/G_0^0=S^3\times S^3/\triangle
$$
o{\`u} $\triangle$ est la diagonale.}
\end{rem}

Soit $\g ,\g_0,\g_2$ les alg{\`e}bres de Lie respectives de $G^0,G_0^0$
et $G_2^0$. Alors on $\g=\g_2\oplus\g_0$ et $\g_2$ est un id{\'e}al
de $\g$ (on a $[\g,\g_2]\subset\g_2$) et est stable sous l'action
adjointe de $G^0$.
On a $\g=\left\{
  \begin{pmatrix}
    R_{\alpha} +L_{\delta} & 0 \\
    0 & R_{\beta} + L_{\delta}
  \end{pmatrix},\alpha,\,\beta,\,\delta\in \im\h\right\}$,
  $\g_0=\left\{\begin{pmatrix}
    R_{\alpha} +L_{\delta} & 0 \\
    0 & R_{\alpha} + L_{\delta}
  \end{pmatrix},\alpha,\,\delta\in \im\h\right\}$, $\g_2=\left\{
  \begin{pmatrix}
    0 & 0 \\
    0 & R_{\gamma}
    \end{pmatrix},\gamma\in \im\h\right\}$.\\

\vspace{0.7mm}
Consid{\'e}rons le groupe $\G$, composante neutre du groupe des
isom{\'e}tries affines de $\R^8$ conservant la nullit{\'e} de $B$, que l'on
repr{\'e}sente comme $G^0\ltimes \R^8$ muni du produit
$$
(G,T)\cdot (G',T')=(GG',GT' + T).
$$
Alors l'alg{\`e}bre de lie $\gt$ de $\G$ s'{\'e}crit $\gt=\g\oplus\R^8=\g_0\oplus
\g_2\oplus\R^8$, le crochet {\'e}tant donn{\'e} par
$$
[(\eta,t),(\eta',t')]=([\eta,\eta'],\eta t' - \eta't)
$$
On a alors les relations suivantes :
$[\g, \R^8]= \R^8$, $[ \R^8, \R^8]=0$, $[\g_0,\g_0]=\g_0$, $[\g_0,\g_2]=\g_2$,
$[\g_2,\g_2]=\g_2$.


\section{Surfaces $\Sigma_V$}

\subsection{Immersions conformes $\Sigma_V$}

\begin{defi}
On dira qu'une surface immerg{\'e}e, $\Sigma$, de $\oct$ est \textbf{une
surface} $\Sigma_V$ si $\forall z\in \Sigma,\ T_z\Sigma\in Q$. En outre {\`a}
$\Sigma$ est associ{\'e}e la fonction $\rho_{\Sigma}$ {\`a} valeurs dans $S^3$ d{\'e}finie
 par
$$
\rho_{\Sigma}(z)=\rho(T_z\Sigma)
\ .$$
En particulier, soit $X\colon \Omega \to \oct$ une immersion
conforme d'un ouvert simplement connexe de $\R^2$ dans $\oct$, alors
on dira que c'est une immersion conforme $\Sigma_V$ si
$$
\forall z=(u,v)\in \Omega,\quad dX=e^f(q\,du + q'dv)
$$
avec $(q,q')\in V$, i.e. $|q|=|q'|=1$ et $B(q,q')=0$.
En outre on dira que $z\in \Sigma$ est \textbf{un point r{\'e}gulier} de $\Sigma$
si $T_z\Sigma\in U$ (i.e. $0<|\im\!(x.\bar{x'})|<\frac{1}{2}$ avec $q=(x,y),\,
q'=(x',y')$, pour une immersion conforme $\Sigma_V$). Dans le cas
contraire, on dira que $\Sigma $ admet \textbf{un  point singulier} en $z$.
On dira alors  que c'est \textbf{un point singulier de type} $P_1$ si
$T_z\Sigma\in G.P_1$, et de \textbf{type} $P_2$ si $T_z\Sigma\in G.P_2$
(i.e. $|\im\!(x.\bar{x'})|=0$ et $|\im\!(x.\bar{x'})|=\frac{1}{2}$ respectivement).
\end{defi}

\begin{defi}
On appellera rel{\`e}vement $\Sigma_V$ une application $U=(F,X)\colon
\Omega \to \G$ telle que $X$ soit une immersion conforme $\Sigma_V$ et
que $\tilde{\rho}\circ F=\rho_X$.
\end{defi}

Le groupe de gauge $C^{\infty}(\Omega,G_0^0)$ agit sur l'ensemble
$\G(\Sigma_V)$ des rel{\`e}vements $\Sigma_V$ : $(F,X)\cdot(K,0)=(FK,X)$.
 L'orbite de $(F,X)$ est l'ensemble des rel{\`e}vements correspondants au m{\^e}me $X$.
  Dans chaque orbite, on peut prendre par exemple
$$
F=\mathcal{R}_{\rho_X}:=
  \begin{pmatrix}
    Id & 0 \\
    0 & R_{\rho_X}
  \end{pmatrix}
$$
alors tout rel{\`e}vement de $X$ est de la forme $(\mathcal{R}_{\rho_X} M,X)$
avec $M\in C^{\infty}(\Omega,G_0^0)$.

\subsubsection{Forme de Maurer-Cartan}

Soit $U=(F,X)=(\mathcal{R}_{\rho} M,X)$ un rel{\`e}vement $\Sigma_V$ alors sa
forme de Maurer-Cartan est donn{\'e}e par $U^{-1}.dU=(F^{-1}.dF,F^{-1}.dX)$ avec
\begin{eqnarray}\label{f-1df}
\negthickspace\bullet \ F^{-1}.dF & = & M^{-1}
  \begin{pmatrix}
    0 & 0 \\
    0 & R_{d\rho .\rho^{-1}}
  \end{pmatrix}M + M^{-1}.dM \hfill\null\\
  & = &  \begin{pmatrix}
    0 & 0 \\
    0 & R_{a(d\rho .\rho^{-1})a^{-1}}
  \end{pmatrix} +
  \begin{pmatrix}
    R_{da.a^{-1}} + L_{c^{-1}dc} & 0 \\
    0 & R_{da.a^{-1}} + L_{c^{-1}dc}
  \end{pmatrix}\nonumber
\end{eqnarray}
 en posant $M=Diag(R_aL_c,R_aL_c)$ ($(a,c)$ n'est d{\'e}fini qu'{\`a} $\pm 1$
 pr{\`e}s mais $\Omega$ est simplement connexe.),\\

\noindent$\bullet \ F^{-1}.dX=e^f(E_1du + E_2dv)$ avec $\rho(E_1,E_2)=1$
 (d'apr{\`e}s $\tilde{\rho}\circ F=\rho_X$ et le th{\'e}or{\`e}me \ref{semidirect}-(iv)).
Ainsi $E_2= {0\choose 1}.E_1=\flb[E].E_1$ d'o{\`u}
$$
F^{-1}.dX=
  \begin{pmatrix}  x \\ y   \end{pmatrix}du + \begin{pmatrix}  -y \\  x
  \end{pmatrix}dv \quad \text{avec}\quad |x|^2 + |y|^2=e^{2f}.
$$
Reciproquement, si $F^{-1}.dX$ est de cette forme, alors $X$ est une
immersion conforme $\Sigma_V$ et d'apr{\`e}s le th{\'e}or{\`e}me \ref{semidirect}-(iv),
 $U=(F,X)$ est un rel{\`e}vement $\Sigma_V$, i.e.
 $\rho_X=\tilde{\rho}\circ F$. D'o{\`u} le th{\'e}or{\`e}me suivant et son
 corollaire :

\begin{theo}
Soit $U=(F,X)\colon \Omega \to \G$ , alors $U\in \G(\Sigma_V)$ (i.e.
$X$ est une immersion conforme $\Sigma_V$ et $\rho_X=\tilde{\rho}\circ
F$) \ssi
$$
F^{-1}.dX=
  \begin{pmatrix}  x \\ y   \end{pmatrix}du + \begin{pmatrix}  -y \\  x
  \end{pmatrix}dv \quad \text{avec}\quad  (x,y)\neq 0.
$$
En outre $z_0$ est un point r{\'e}gulier \ssi
$0<|\im\!(x.\bar{x'})|<\frac{1}{2}(|x|^2 + |y|^2)$, un point singulier de
type $P_1$ ou $P_2$ \ssi $|\im\!(x.y)|=0$ ou $|\im\!(x.y)|=\frac{1}{2}(|x|^2 +
 |y|^2)$ respectivement.
\end{theo}

\begin{cor}
Soit $\alpha\in T^*\Omega\otimes\gt$, alors $\alpha $ est la forme
de Maurer-Cartan d'un {\'e}l{\'e}ment $U\in\G(\Sigma_V)$ \ssi :
\begin{description}
  \item[(i)] $d\alpha + \alpha \wedge\alpha=0$
  \item[(ii)] si on pose $\alpha=(\eta,t)$ alors $t=\begin{pmatrix}  x \\ y
 \end{pmatrix}du + \begin{pmatrix}  -y \\  x
  \end{pmatrix}dv$ \ avec  $(x,y)\neq 0$.
\end{description}
Dans ce cas,  suivant les valeurs de $|\im\!(x.y)|$, on peut
connaitre le type du point $z\in\Omega$ pour l'immersion $X$.
\end{cor}

Posons
\begin{eqnarray*}
\g_{-1} & = & \left\{\begin{pmatrix} x \\   y  \end{pmatrix} -i \begin{pmatrix}
  -y \\  x  \end{pmatrix},\, (x,y)\in \oct\right\}\subset \oct\otimes\C \\
\g_{1} & = & \left\{\begin{pmatrix} x \\   y  \end{pmatrix} +i \begin{pmatrix}
 -y \\  x  \end{pmatrix},\, (x,y)\in \oct\right\}\subset \oct\otimes\C .
\end{eqnarray*}
Ce sont des sous-espaces complexes de $\oct\otimes\C$ : ce sont les
sous-espaces propres associ{\'e}s aux valeurs propres $i$ et $-i$
respectivement pour l'endomorphisme de $\oct\otimes\C$ : $L_{\flb[E]}$.
En particulier $\oct\otimes\C=\g_{-1}\oplus\g_1$ et
$\g_{-1}=\overline{\g_1}$.\\
Les actions respectives de $(\R_{+}^*)^2$ et $SO(2)$ stabilisent
l'ensemble des couples
$(\genfrac(){0pt}{1pt}{x}{y} , \genfrac(){0pt}{1pt}{-y}{x})$:
\begin{eqnarray*}
(\alpha,\beta)\cdot \left(\begin{pmatrix} x \\   y  \end{pmatrix},\begin{pmatrix}
 -y \\  x  \end{pmatrix}\right) & = & \left(\begin{pmatrix} \alpha x \\  \beta y
   \end{pmatrix},\begin{pmatrix} -\beta y \\  \alpha x
   \end{pmatrix}\right)\\
\left(\begin{pmatrix} x \\   y  \end{pmatrix} ,\begin{pmatrix}
 -y \\  x  \end{pmatrix}\right)\cdot R_{\theta} & = & \left(R_{\theta}\begin{pmatrix}
  x \\   y  \end{pmatrix}, R_{\theta}\begin{pmatrix}
  -y \\   x  \end{pmatrix}\right).
\end{eqnarray*}
Ainsi, on peut faire agir ces deux groupes sur $\g_{-1}$ et $\g_1$
respectivement et donc sur $\oct\otimes\C$. En particulier, on a
\begin{eqnarray*}
R_{\theta}\left(\begin{pmatrix}
  x \\   y  \end{pmatrix} -i\begin{pmatrix}
  -y \\   x  \end{pmatrix}\right) & = & e^{i\theta}\left(\begin{pmatrix}
  x \\   y  \end{pmatrix} -i\begin{pmatrix}
  -y \\   x  \end{pmatrix}\right)\\
R_{\theta}\left(\begin{pmatrix}
  x \\   y  \end{pmatrix} +i\begin{pmatrix}
  -y \\   x  \end{pmatrix}\right) & = & e^{-i\theta}\left(\begin{pmatrix}
  x \\   y  \end{pmatrix} +i\begin{pmatrix}
  -y \\   x  \end{pmatrix}\right).
\end{eqnarray*}
Posons
\begin{eqnarray*}
\g_{-1}^* & = & \left\{\begin{pmatrix} x \\   y  \end{pmatrix} -i \begin{pmatrix}
  -y \\  x  \end{pmatrix} / 0<|\im\!(x.y)|<\frac{1}{2}(|x|^2 + |y|^2)\right\} \\
\g_{1}^* & = & \left\{\begin{pmatrix} x \\   y  \end{pmatrix} +i \begin{pmatrix}
 -y \\  x   \end{pmatrix} /  0<|\im\!(x.y)|<\frac{1}{2}(|x|^2 +  |y|^2)\right\}.
\end{eqnarray*}
Ce sont des ouverts, stables par homoth{\'e}tie complexe, de $\g_{-1}$
et $\g_1$ respectivement. Posons aussi $F_1=\{q-iL_{\flb[E]}
q/\,\im\!(x.y)=0\}$ et
$F_2=\{q-iL_{\flb[E]}q/\,|\im\!(x.y)|$ $=$ $\frac{1}{2}(|x|^2 +  |y|^2)\}$.\\
Soit $\alpha\in T^*\Omega\otimes \gt$, $\alpha=(\eta,t)$, tel que
$d\alpha + \alpha\wedge \alpha =0$ et {\'e}crivons $t=\alpha_{-1} +
\alpha_1$ la d{\'e}composition de $t$ suivant $\g_{-1}\oplus\g_1$. Alors
on a $\alpha_1=\overline{\alpha_{-1}}$ car $t$ est r{\'e}elle. Alors
$\alpha$ est la forme de Maurer-Cartan d'un {\'e}l{\'e}ment de $\G(\Sigma_V)$ \ssi
$\alpha_{-1}=\alpha_{-1}({\partial \over \partial z})dz$ (et donc
$\alpha_1=\alpha_1(\frac{\partial}{\partial \bar{z}})d\bar{z}$).
Ceci permet de r{\'e}{\'e}crire le corrollaire sous la forme:

\begin{theo}
Soit $\alpha\in T^*\Omega\otimes \gt$, tel que $d\alpha + \alpha\wedge
 \alpha =0$. Alors $\alpha$ correspond {\`a} un {\'e}l{\'e}ment $\G(\Sigma_V)$
 \ssi $\alpha_{-1}(\frac{\partial}{\partial \bar{z}})=0$ et $\alpha_{-1}
 (\frac{\partial}{\partial z})$ ne s'annule pas . Dans ce cas $X$ a
 un point r{\'e}gulier en $z_0$ \ssi $\alpha_{-1} (\frac{\partial}{\partial z})
 (z_0)\in \g_{-1}^*$, un point singulier de type $P_1$ ou $P_2$ \ssi
$\alpha_{-1} (\frac{\partial}{\partial z})(z_0)\in F_1$ ou $F_2$
respectivement.
\end{theo}

\begin{rem}{\em
\begin{enumerate}
  \item L'action du groupe de gauge $C^{\infty}(\Omega,G_0^0)$ sur
  $\G(\Sigma_V)$ induit une action sur les formes de Maurer-Cartan:
  $$(\eta,t)\mapsto (K\eta K^{-1}-dK.K^{-1},\,K.t)$$
  \item $\oct\otimes \C=\h^2\otimes\C$ est un $\h$-espace vectoriel
  ({\`a} gauche ou {\`a} droite au choix). Il en est de m{\^e}me de $\g_{-1}$
  et de $\g_1$:
\begin{eqnarray*}
\g_{-1} & = & \{x.\epsilon +
y.(L_{\flb[E]}.\epsilon),\,(x,y)\in\h^2\}=
\h .\epsilon\oplus\h .(L_{\flb[E]}.\epsilon) \\
\g_{1} & = &  \{x.\bar{\epsilon} +
y.(L_{\flb[E]}.\bar{\epsilon}),\,(x,y)\in\h^2\}=
\h .\bar{\epsilon}\oplus\h .(L_{\flb[E]}.\bar{\epsilon})
\end{eqnarray*}
o{\`u} $\epsilon=\displaystyle\frac{1}{2}  \begin{pmatrix} 1 \\   -i  \end{pmatrix}.$
  \item On a aussi $\oct\otimes\C=(\h\otimes\C)^2$. Ainsi
\begin{eqnarray*}
\g_{-1} & = & \{(x+iy).\epsilon\,
/(x+iy)\in\h\otimes\C\}=(\h\otimes\C).\epsilon \\
\g_{1} & = &  \{(x+iy).\bar{\epsilon}\,/(x+iy)\in\h\otimes\C\}=(\h\otimes\C).
\bar{\epsilon}.
\end{eqnarray*}
\end{enumerate}}
\end{rem}

\subsection{D{\'e}composition de l'alg{\`e}bre de Lie.}

Consid{\'e}rons l'automorphisme int{\'e}rieur $\tau$ de $G$ d{\'e}fini par
$(-L_{\flb[E]},0)$ :
$$
\tau(G,T)=(-L_{\flb[E]},0)(G,T)(-L_{\flb[E]},0)^{-1}=(-L_{\flb[E]}G
L_{\flb[E]},-L_{\flb[E]}T)
$$
$\tau$ induit un automorphisme de l'alg{\`e}bre de Lie $\gt$
($=Ad_{(-L_{\flb[E]},0)}$) qui v{\'e}rifie $\tau^4=Id$, donc
$\tau$ est diagonalisable dans $\gt^{\C}=\gt\otimes \C$ et ses
valeurs propres sont les $i^k,\, 0\leq k\leq 3$. On notera
$\gt_{k}^{\C}$ les espaces propres. On a alors :
\begin{itemize}
  \item[$\bullet$] $\gt_{-1}^{\C}=\g_{-1}$
  \item[$\bullet$] $\gt_{1}^{\C}=\g_{1}$
  \item[$\bullet$]$\gt_0^{\C}=\g_0^{\C}=\g_0\otimes \C$
  \item[$\bullet$]$\gt_{2}^{\C}=\gt_{2}\otimes \C $\,, o{\`u}
  $\gt_2=\{Diag(-R_{\gamma},\,R_{\gamma}),\,\gamma\in \im\h\}$ n'est
  pas une alg{\`e}bre de Lie.
\end{itemize}
Comme $\tau$ est un automorphisme on
a $[\gt_k^{\C},\gt_l^{\C}]\subset \gt_{(k+l)\ \text{mod }4}^{\C}$ \ et
de plus on a aussi $[\gt_{\pm 1}^{\C},\gt_{\pm 1}^{\C}]=0$.\\
On notera $[.]_k\colon\gt^{\C} \to \gt_k^{\C}$ la projection sur
$\gt_k^{\C}$ suivant la d{\'e}composition $\gt^{\C}=\gt_{-1}^{\C}\oplus\gt_0^{\C}
\oplus\gt_1^{\C}\oplus\gt_2^{\C}$, et $\alpha_k=[\alpha]_k$. Alors
on a
$$
\alpha=\alpha_{-1} + \alpha_0 +\alpha_1 +\alpha_2
$$
En substituant cette expression de $\alpha$ dans l'{\'e}quation $d\alpha
+ \alpha\wedge\alpha=0$ on obtient en projetant le r{\'e}sultat sur
chaque espace propre :
\begin{equation}\label{alpha-k}
\left\{\begin{array}{ccc}
d\alpha_{-1} + [\alpha_{-1}\wedge\alpha_0] + [\alpha_1\wedge\alpha_2] & =
 & 0\\
d\alpha_0 + \frac{1}{2}[\alpha_0\wedge\alpha_0] + \frac{1}{2}[\alpha_2\wedge
\alpha_2] & = & 0\\
d\alpha_{1} + [\alpha_{1}\wedge\alpha_0] + [\alpha_{-1}\wedge\alpha_2] & =
 & 0\\
d\alpha_2 +[\alpha_0\wedge\alpha_2] & = & 0
\end{array}\right.
\end{equation}
Posons $\alpha_k'=\alpha_k(\frac{\partial}{\partial z})dz$,
$\alpha_k''=\alpha(\frac{\partial}{\partial \bar{z}})d\bar{z}$. On
a vu que $\alpha$ est la forme de Maurer-Cartan d'un {\'e}l{\'e}ment de
$\G(\Sigma_V)$ \ssi $\alpha_{-1}=\alpha_{-1}'$ et
$\alpha_1=\alpha_1''$. Ainsi $\alpha =\alpha_2' + \alpha_{-1}' + \alpha_0 +
\alpha_1''+ \alpha_2''$. Les {\'e}quations (\ref{alpha-k}) s'{\'e}crivent
alors :
\begin{equation}
\left\{\begin{array}{ccc}
d\alpha_{-1}' + [\alpha_{-1}'\wedge\alpha_0] + [\alpha_1''\wedge\alpha_2'] & =
 & 0\\
d\alpha_0 + \frac{1}{2}[\alpha_0\wedge\alpha_0] + \frac{1}{2}[\alpha_2'\wedge
\alpha_2''] & = & 0\\
d\alpha_{1}'' + [\alpha_{1}''\wedge\alpha_0] + [\alpha_{-1}'\wedge\alpha_2''] & =
 & 0\\
d\alpha_2 +[\alpha_0\wedge\alpha_2] & = & 0
\end{array}\right.
\end{equation}
Posons pour $\lambda\in\C^*$,
\begin{eqnarray*}
\alpha_{\lambda}& = & \lambda^{-2}\alpha_2' + \lambda^{-1}\alpha_{-1}' +
\alpha_0 + \lambda\alpha_1'' + \lambda^2\alpha_2'' \\
 & = & \lambda^{-1}\alpha_{-1}' + \alpha_0 + \lambda\alpha_1'' +
\frac{\lambda^2 + \lambda^{-2}}{2}\,\alpha_2 +
\frac{\lambda^2 - \lambda^{-2}}{2i}(\ast\alpha_2).
\end{eqnarray*}
Alors on a
\begin{eqnarray*}
d\alpha_{\lambda} + \alpha_{\lambda}\wedge\alpha_{\lambda} & = &
\lambda^{-1}(d\alpha_{-1}' + [\alpha_{-1}'\wedge\alpha_0] +
[\alpha_1''\wedge\alpha_2']) \\
 &   &  +(d\alpha_0 + \frac{1}{2}[\alpha_0\wedge\alpha_0] + \frac{1}{2}[\alpha_2'\wedge
\alpha_2''])\\
 &  & + \lambda(d\alpha_{1}'' + [\alpha_{1}''\wedge\alpha_0] +
 [\alpha_{-1}'\wedge\alpha_2'']) \\
 &  & +\frac{\lambda^2 + \lambda^{-2}}{2}(d\alpha_2 +
 [\alpha_0\wedge\alpha_2])\\
 &  & +\frac{\lambda^2 - \lambda^{-2}}{2i}(d(\ast\alpha_2) +
 [\alpha_0\wedge(\ast\alpha_2)])\\
 & = & \frac{\lambda^2 - \lambda^{-2}}{2i}(d(\ast\alpha_2) +
 [\alpha_0\wedge(\ast\alpha_2)])
\end{eqnarray*}
On peut maintenant d{\'e}montrer le th{\'e}or{\`e}me suivant:

\begin{theo}\label{integrable}
On suppose $\Omega$ simplement connexe. Soit $\alpha \in
T^*\Omega\otimes\gt$. Alors
\begin{itemize}
  \item[$\bullet$] $\alpha$ est la forme de Maurer-Cartan d'un
  {\'e}l{\'e}ment de $\G(\Sigma_V)$ \ssi $d\alpha + \alpha \wedge\alpha=0$,
  $\alpha_{-1}''=\alpha_1'=0$ et $\alpha_{-1}',\alpha_1''$ ne
  s'annule pas.
  \item[$\bullet$] Dans ce cas, $\alpha$ correspond {\`a} une immersion
  conforme  $\Sigma_V$ telle que $\rho_X$ est harmonique \ssi la forme
  de Maurer-Cartan prolong{\'e}e $\alpha_{\lambda}=\lambda^{-2}\alpha_2' +
  \lambda^{-1}\alpha_{-1}' +\alpha_0 + \lambda\alpha_1'' +
  \lambda^2\alpha_2''$ v{\'e}rifie
\begin{equation}\label{courburnul}
d\alpha_{\lambda} + \alpha_{\lambda}\wedge\alpha_{\lambda}=0,\
\forall \lambda\in \C^*.
\end{equation}
\end{itemize}
\end{theo}

\noindent{\em D{\'e}monstration} --- On a d{\'e}j{\`a} vu le premier point.
Quand au second, il s'agit d'apr{\`e}s le calcul pr{\'e}c{\'e}dent de montrer que:
$\rho_X $  est harmonique \ssi \ $d(\ast\alpha_2) +
[\alpha_0\wedge(\ast\alpha_2)]=0 $.
 Or on a
$$\alpha_0 +\alpha_2=F^{-1}.dF=
  \begin{pmatrix}
    0 & 0 \\
    0 & R_{a(d\rho.\rho^{-1})a^{-1}}
  \end{pmatrix} + M^{-1}.dM$$
 d'o{\`u}
$$\alpha_2=\frac{1}{2}
  \begin{pmatrix}
    -R_{a(d\rho.\rho^{-1}).a^{-1}} & 0 \\
     0 & R_{a(d\rho.\rho^{-1}).a^{-1}}
  \end{pmatrix}=M^{-1}
  \begin{pmatrix}
    -R_{\gamma/2} & 0 \\
     0 & R_{\gamma/2}
  \end{pmatrix}M$$
 avec $\gamma=d\rho.\rho^{-1}$, et
 $$\alpha_0=M^{-1}dM +
  M^{-1}.\begin{pmatrix}
     R_{\gamma/2} & 0 \\
     0 & R_{\gamma/2}
    \end{pmatrix}M$$
Posons $\beta= \begin{pmatrix}
    -R_{\gamma/2} & 0 \\
     0 & R_{\gamma/2}
  \end{pmatrix}$. Alors $\rho$ est harmonique \ssi $d(\ast\gamma)=
  (\triangle\rho + |d\rho|^2 \rho)\rho^{-1}du\wedge dv=0$
  \ssi $d(\ast\beta)=0$. Or
$$d(\ast\beta)=d(M(\ast\alpha_2)M^{-1})=M(\,d(\ast\alpha_2) + \,
[(M^{-1}.dM)\wedge(\ast\alpha_2)]\ )M^{-1}$$
De plus on a
$[(M^{-1}.dM)\wedge(\ast\alpha_2)]=[\alpha_0\wedge(\ast\alpha_2)]$
puisque
$$[(M^{-1}\begin{pmatrix}
     R_{\gamma/2} & 0 \\
     0 & R_{\gamma/2}
    \end{pmatrix}M)\wedge (\ast\alpha_2)]=\frac{1}{4}M^{-1}
  \begin{pmatrix}
    -R_{[\gamma\wedge(\ast\gamma)]} & 0 \\
    0 & R_{[\gamma\wedge(\ast\gamma)]}
  \end{pmatrix}M=0$$
  car $[\gamma\wedge(\ast\gamma)]=0$.
Finalement on a
$$
d(\ast\beta)=M(\,d(\ast\alpha_2) + [\alpha_0\wedge(\ast\alpha_2)]\ )M^{-1}
$$
ce qui ach{\`e}ve la d{\'e}monstration du th{\'e}or{\`e}me.\hfill $\blacksquare$ \\

\begin{rem}{\em
Chaque point de la surface sera r{\'e}gulier, respectivement singulier de type
 $P_1$, ou $P_2$ \ssi $\alpha_{-1}(\frac{\partial}{\partial z})
\in \g_{-1}^*$, $F_1$ ou $F_2$ respectivement. En outre il suffit
que (\ref{courburnul}) soit vrai pour tout $\lambda \in S^1$ pour
que $\rho_X$ soit harmonique.}
\end{rem}
\begin{cor}
Supposons que $\Omega$ soit simplement connexe. Soit $\alpha\in
T^*\Omega\otimes\gt$ la forme de Maurer-Cartan associ{\'e}e {\`a} une
immersion conforme $\Sigma_V$ dont le $\rho_X$ est harmonique, et
$z_0\in \Omega$. Alors pour tout $\lambda\in S^1$, il existe un unique
rel{\`e}vement $\Sigma_V$, $U_{\lambda}\in C^{\infty}(\Omega,\G)$ tel que
$$
dU_{\lambda}=U_{\lambda}\alpha_{\lambda}\quad \text{et} \quad
U_{\lambda}(z_0)=\mathbf{1}.
$$
Ainsi il existe un famille $(X_{\lambda})_{\lambda\in S^1}$
d'immersions $\Sigma_V$ dont le $\rho_X$ est harmonique donn{\'e}e par
$U_{\lambda}=(F_{\lambda},X_{\lambda})$, tel que $X=X_1$ (en supposant que
$X(z_0)=0$).
En outre si $X$ admet un point r{\'e}gulier (resp. singulier de type
$P_1$ ou $P_2$) en $z\in\Omega$, il en est de m{\^e}me pour
$X_{\lambda}$ pour tout $\lambda\in S^1$. Autrement dit le type d'un
point $z\in \Omega $ est le m{\^e}me pour toutes les immersions
$X_{\lambda}$.
\end{cor}

\noindent{\em D{\'e}monstration} ---  Il suffit d'appliquer le th{\'e}or{\`e}me pr{\'e}c{\'e}dent
et de remarquer que pour $\lambda\in S^1,\, \alpha_{\lambda}$ est {\`a} valeurs
 dans $\gt$, et qu'on a  $(\alpha_{\lambda})_{-1}'=\lambda^{-1}\alpha_{-1}'$.
\hfill $\blacksquare$ \\

\subsection{Equations associ{\'e}es (lin{\'e}aire et non lin{\'e}aire)}

Soit $X$ une immersion $\Sigma_V$ sur $\Omega$ simplement connexe.
Posons $(E_1,E_2)=e^f \mathcal{R}_{\rho}^{-1}(q,q')$. Alors $E_2=\flb[E].E_1$
i.e. en posant $E_1={x \choose y}\neq0$ on a $E_2={ -y \choose x}$.
Alors {\'e}crivons que $dX=\mathcal{R}_{\rho}(E_1 du + E_2 dv)$ est ferm{\'e}e, on
obtient
\begin{equation}\label{du-dv}
0=\frac{\partial E_1}{\partial v}-\frac{\partial E_2}{\partial u} +
  \begin{pmatrix}
    0 & 0 \\
    0 & R_{\gamma_v}
  \end{pmatrix}E_1 -\begin{pmatrix}
    0 & 0 \\
    0 & R_{\gamma_u}
  \end{pmatrix}E_2
\end{equation}
o{\`u} $\gamma=\gamma_u du + \gamma_v dv=d\rho.\rho^{-1}$, ce qui
s'{\'e}crit encore
\begin{equation}\label{du-dv'}
 \left\{ \begin{array}{l}
    \displaystyle \frac{\partial x}{\partial v} + \frac{\partial y}
    {\partial u}=0\\
\displaystyle \frac{\partial y}{\partial u}- \frac{\partial x}{\partial v} +
y.\gamma_v - x.\gamma_u=0
  \end{array}\right.
\end{equation}
On va essayer de retrouver cette {\'e}quation en utilisant le rel{\`e}vement
$U=(\mathcal{R}_\rho,X)$. Alors on a
$$
\alpha =U^{-1}.dU=\left(\begin{pmatrix}
    0 & 0 \\
    0 & R_{\gamma}
  \end{pmatrix},E_1 du +E_2 dv\right)
  .$$
Posons $E=\frac{1}{2}(E_1-iE_2)=\frac{1}{2}({x \choose y}-i{-y
\choose x})=x.\epsilon + y.(L_{\flb[E]}.\epsilon)=(x+iy).\epsilon$.
On a alors :
$$
\alpha_{-1}'=Edz,\ \alpha_1''=\bar{E}.d\bar{z},\ \alpha_0=Diag\left(\frac{1}{2}
R_{\gamma},\frac{1}{2}R_{\gamma}\right),\ \alpha_2=Diag\left(-\frac{1}{2}
R_{\gamma},\frac{1}{2}R_{\gamma}\right).$$
Projetons l'{\'e}quation  $d\alpha + \alpha\wedge\alpha=0$ sur
$\g_{-1}$, on obtient:
\begin{equation}\label{dedzbar}
\frac{\partial E}{\partial \bar{z}} + \frac{1}{2}
 \begin{pmatrix}
   R_{\gamma(\frac{\partial}{\partial \bar{z}})} & 0 \\
    0 & R_{\gamma(\frac{\partial}{\partial \bar{z}})}
  \end{pmatrix}.E + \frac{1}{2}\begin{pmatrix}
   R_{\gamma(\frac{\partial}{\partial z})} & 0 \\
    0 & -R_{\gamma(\frac{\partial}{\partial z})}
  \end{pmatrix}.\bar{E}=0 \ .
\end{equation}
On v{\'e}rifie facilement que cette {\'e}quation est {\'e}quivalente {\`a}
(\ref{du-dv}). Nous allons maintenant utiliser le th{\'e}or{\`e}me
\ref{alphabeta} pour r{\'e}{\'e}crire cette {\'e}quation.\\
Soit $g:\Omega\to G$, $\alpha,\beta: \Omega\to (\R_+^{\ast})^2$, et
$\theta:\Omega\to \R$ tel que
$$
(E_1,E_2)=(g\cdot(\alpha,\beta)\cdot(q_0,q_0'))\cdot R_{\theta}
$$
o{\`u} $(q_0,q_0')=\left(\begin{pmatrix}
    \frac{1}{\sqrt{2}}  \\ \frac{I}{\sqrt{2}}
      \end{pmatrix},
  \begin{pmatrix}
   \frac{-I}{\sqrt{2}}  \\  \frac{1}{\sqrt{2}}  \end{pmatrix}\right)$ par
exemple, et o{\`u} on {\'e}crit $\h=\R\oplus \R I\oplus\R J\oplus \R K$
pour ne pas confondre le $i$ des complexes provenant de la complexification
avec celui de $\h$. Comme $\rho(q_0,q_0')=\rho(E_1,E_2)=1$, on
a $\tilde{\rho}(g)=1$ i.e. $g=Diag(R_aL_c,R_aL_c)$. On peut aussi
{\'e}crire que
$E=\frac{1}{2}(E_1-iE_2)=e^{i\theta}g.((\alpha x_0+i\beta y_0).\epsilon)$
alors (\ref{dedzbar}) s'{\'e}crit
$$
  \begin{array}{c}
e^{i\theta}\displaystyle\left[\,i\frac{\partial\theta}{\partial\bar{z}}
\,g\cdot((\alpha x_0 + i\beta y_0)\cdot\epsilon) +
\frac{\partial g}{\partial \bar{z}}
\cdot((\alpha x_0 +i\beta y_0)\cdot\epsilon)\right. + \\
\displaystyle \left. g\cdot\left(\left(\frac{\partial \alpha}{\partial\bar{z}}
x_0 + i\frac{\partial\beta}{\partial \bar{z}}y_0\right)\cdot\epsilon\right)
\,\right] + \displaystyle\frac{1}{2}e^{i\theta}
  \begin{pmatrix}
    R_{\gamma(\frac{\partial}{\partial\bar{z}})} & 0 \\
    0 & R_{\gamma(\frac{\partial}{\partial\bar{z}})}
  \end{pmatrix}\cdot g \cdot((\alpha x_0 + i\beta y_0)\cdot\epsilon)\\
+ \displaystyle\frac{1}{2}e^{-i\theta}
  \begin{pmatrix}
    R_{\gamma(\frac{\partial}{\partial z})} & 0 \\
    0 & -R_{\gamma(\frac{\partial}{\partial z})}
  \end{pmatrix}\cdot g\cdot((\alpha x_0 - i\beta y_0)\cdot\bar{\epsilon})=0
 \end{array}
 $$
 ce qui s'{\'e}crit encore, tout calcul fait et en factorisant par
 $e^{i\theta}.g$ :
\begin{equation}\label{thetalphabeta}
\begin{array}{c}
i{\displaystyle\frac{\partial\theta}{\partial\bar{z}}}A +
\delta(\frac{\partial}{\partial\bar{z}})A + A\,\tilde{\alpha}(\frac{\partial}
{\partial\bar{z}})+ {\displaystyle\frac{\partial A}{\partial\bar{z}}} \\
+ {\displaystyle\frac{1}{2}A}
\left(a\gamma(\frac{\partial}{\partial \bar{z}})a^{-1}\right) +
{\displaystyle \frac{e^{-2i\theta}}{2}\bar{A}}
\left(a\gamma(\frac{\partial}{\partial z})a^{-1}\right)=0.
\end{array}
\end{equation}
o{\`u} on a pos{\'e} $A=\alpha x_0 + i\beta y_0=\frac{1}{\sqrt{2}}(\alpha +
iI\beta)$, $\tilde{\alpha}=da.a^{-1}$, $\delta=c^{-1}.dc$\,.
 Les inconnues $\tilde{\alpha}$, $\delta$ et $\gamma'=a\gamma a^{-1}$
sont les param{\'e}tres qui interviennent dans la forme de Maurer-Cartan,
$F^{-1}.dF$, du rel{\`e}vement $F$ de $\rho_X$, donn{\'e}e par (\ref{f-1df}).
Ainsi, on peut consid{\'e}rer que l'on construit $\rho_X$ {\`a} partir de la
repr{\'e}sentation de Weierstrass pour l'espace sym{\'e}trique $S^3=G^0/G_0^0$
(cf. \cite{DPW}), alors cela nous donne la forme de Maurer-Cartan, $F^{-1}.dF$, et on
peut donc consid{\'e}rer  $\tilde{\alpha}$, $\delta$ et $\gamma'$ comme
des param{\`e}tres, les inconnues {\'e}tant alors $\theta$ et $A$.
Cependant, il y a alors un probl{\`e}me de compatibilit{\'e} puisque
$A\in (\R\oplus\R I)\otimes\C$.


\section{Groupe de lacets}

\subsection{D{\'e}finitions et notations}

\begin{defi}
Soit $G$ un groupe de Lie, on appellera groupe de lacets sur $G$, le
groupe $C^{\infty}(S^1,G)$ que l'on notera $\Lambda G$ (cf. \emph{\cite{PS}}).
\end{defi}
Dans notre cas, les groupes consid{\'e}r{\'e}s sont
$\G,\,G^0,\,G_0^0,\,G_2^0$. On d{\'e}finit les groupes suivants:
\begin{eqnarray*}
\Lambda\G_{\tau} & = & \{[\lambda \mapsto U_{\lambda}]\in\Lambda\G
 /U_{i\lambda}=\tau(U_{\lambda})\}\\
\Lambda\G_{\tau}^{\C} & = & \{[\lambda \mapsto U_{\lambda}]\in
 \Lambda\G^{\C} /U_{i\lambda}=\tau(U_{\lambda})\}\\
\Lambda_{\ast}^{-}\G_{\tau}^{\C} & = & \{[\lambda \mapsto U_{\lambda}]\in
 \Lambda\G_{\tau}^{\C} / U_{\lambda}\text{ se prolonge en une fonction
 holomorphe sur}\\
  & & \text{le compl{\'e}mentaire du  disque unit{\'e} et }U_{\infty}=1 \}\\
\Lambda^{+}\G_{\tau}^{\C} & = & \{[\lambda \mapsto U_{\lambda}]\in
 \Lambda\G_{\tau}^{\C} / U_{\lambda}\text{ se prolonge en une
 fonction holomorphe sur}\\
 & & \text{le disque unit{\'e} }\}\\
\Lambda_{\B}^{+}\G_{\tau}^{\C} & = & \{[\lambda \mapsto U_{\lambda}]\in
 \Lambda\G_{\tau}^{\C} / U_{\lambda}\text{ se prolonge en une
 fonction holomorphe sur}\\
 & &  \text{le disque unit{\'e} et }U_0\in (\B,0)\}
\end{eqnarray*}
o{\`u} $\B$ est un sous-groupe de $G^0$. De mani{\`e}re analogue, on
d{\'e}finit les alg{\`e}bres de Lie suivantes:

\begin{eqnarray*}
\Lambda\gt_{\tau}^{\C} & = & \{[\lambda \mapsto \alpha_{\lambda}]\in
\Lambda\gt^{\C}/\alpha_{i\lambda}=\tau(\alpha_{\lambda})\}\\
\Lambda\gt_{\tau} & = & \{[\lambda \mapsto \alpha_{\lambda}]\in
\Lambda\gt_{\tau}^{\C}/
\alpha_{\lambda}\in \gt ,\ \forall \lambda \in S^1\}\\
\Lambda_{\ast}^{-}\gt_{\tau}^{\C} & = & \{[\lambda \mapsto \alpha_{\lambda}]\in
 \Lambda\gt_{\tau}^{\C} / \alpha_{\lambda}\text{ se prolonge en une fonction
 holomorphe sur}\\
  & & \text{le compl{\'e}mentaire du  disque unit{\'e} et }\alpha_{\infty}=0 \}\\
\Lambda^{+}\gt_{\tau}^{\C} & = & \{[\lambda \mapsto \alpha_{\lambda}]\in
 \Lambda\gt_{\tau}^{\C} / \alpha_{\lambda}\text{ se prolonge en une
 fonction holomorphe sur}\\
 & & \text{le disque unit{\'e} }\}\\
\Lambda_{\mathfrak{b}}^{+}\gt_{\tau}^{\C} & = & \{[\lambda \mapsto
\alpha_{\lambda}]\in
 \Lambda\gt_{\tau}^{\C} / \alpha_{\lambda}\text{ se prolonge en une
 fonction holomorphe sur}\\
 & &  \text{le disque unit{\'e} et }\alpha_0\in (\mathfrak{b},0)\}.
\end{eqnarray*}
o{\`u} $\mathfrak{b}$ est une sous-alg{\`e}bre de Lie de $\g_0$.\\
On voit que la condition $\alpha _{i\lambda}=\tau(\alpha_{\lambda}),
\ \forall \lambda \in S^1$ est {\'e}quivalente {\`a}
$\hat{\alpha}_k \in \gt_{k\ \text{mod}\ 4}^{\C}$ \ avec
$\alpha_{\lambda}=\sum_{k\in\mathbb{Z}}\hat{\alpha}_k\lambda^k$.
En outre, on a
$\Lambda\gt_{\tau}^{\C}=\Lambda_{\ast}^{-}\gt_{\tau}^{\C}\oplus\Lambda^{+}
\gt_{\tau}^{\C}$, ce qui permet de d{\'e}finir une projection
:$[.]_{\Lambda_{\ast}^{-}\gt_{\tau}^{\C}}\colon
\Lambda\gt_{\tau}^{\C}\to \Lambda_{\ast}^{-}\gt_{\tau}^{\C}$. On
peut alors r{\'e}ecrire le r{\'e}sultat de la section pr{\'e}c{\'e}dente :

\begin{cor}\label{lacet}
Soit $\alpha$ une 1-forme sur $\Omega$ {\`a} valeurs dans $\gt$ qui
donne lieu {\`a} une immersion $\Sigma_V$ dont le $\rho_X$ est harmonique.
Il lui correspond alors une 1-forme {\`a} valeurs dans
$\Lambda\gt_{\tau}$, $\alpha_{\lambda}$ , qui v{\'e}rifie l'{\'e}quation de
courbure nulle (\ref{courburnul}), et telle que
$$\begin{array}{lcl}
\displaystyle
\left[\alpha_{\lambda}\left(\frac{\partial}{\partial z}\right)
\right]_{\Lambda_{\ast}^{-}
\gt_{\tau}^{\C}} & = & \displaystyle\lambda^{-2}\hat{\alpha}_{-2}\left(\frac
{\partial}{\partial z}\right) +
\lambda^{-1}\hat{\alpha}_{-1}\left(\frac{\partial}{\partial z}\right)\\
\displaystyle\left[\alpha_{\lambda}\left({\partial \over \partial \bar{z}}
\right)\right]_{\Lambda_{\ast}^{-}\gt_{\tau}^{\C}} & = & 0 \\
\end{array}
$$
et
$$
\displaystyle\hat{\alpha}_{-1}\left(\frac{\partial}{\partial z}\right)
 \neq 0.
$$
R{\'e}ciproquement, {\`a} toute 1-forme $\alpha_{\lambda}\in\Lambda\gt_{\tau}$
v{\'e}rifiant ces conditions correspond la 1-forme $\alpha=\alpha_1$ qui
donne lieu {\`a} une immersion $\Sigma_V$ dont le $\rho_X$ est
harmonique. En outre, il existe une unique fonction
$U_{\lambda}\colon \Omega\to \Lambda\G_{\tau}$ tel que
$dU_{\lambda}=U_{\lambda}\alpha_{\lambda}$ et $U_{\lambda}(z_0)=\mathbf{1}$.
 $U_{\lambda}$ sera appel{\'e}e une extention $\Sigma_V$ de $U=U_1$.
\end{cor}

\noindent{\em D{\'e}monstration} --- C'est une cons{\'e}quence imm{\'e}diate du
th{\'e}or{\`e}me \ref{integrable} et du fait que, comme $\alpha$ est
r{\'e}elle, on a $\hat{\alpha}_k=\overline{\hat{\alpha}_{-k}}$.\hfill
$\blacksquare$ \\

\subsection{Th{\'e}or{\`e}mes de d{\'e}composition de groupe}

Ecrivons les d{\'e}compositions d'Iwasawa des diff{\'e}rents groupes que l'on
a rencontr{\'e}s:
\begin{eqnarray*}
\{R_a,\, a\in S^3\}^{\C} & = & \{R_a,\, a\in S^3\}.\B_R\\
\{L_c,\, c\in S^3\}^{\C} & = & \{L_c,\, c\in S^3\}.\B_L .
\end{eqnarray*}
On remarque que $\B_R$ et $\B_L$ commutent. Ensuite, on a
$$
SO(4)^{\C}=\{R_aL_c,\, a,c\in S^3\}^{\C}=SO(4).(\B_R\B_L).
$$
Alors on en d{\'e}duit ${G^0}^{\C}=G^0.\B$ et ${G_0^0}^{\C}=G_0^0.\B_0$
avec
\begin{eqnarray*}
\B & = & \left\{
  \begin{pmatrix}
    AC & 0 \\
    0 & BC
  \end{pmatrix},\, A,B\in \B_R ,\, C\in\B_L\right\}\\
\B_0 & = & \left\{
  \begin{pmatrix}
    AC & 0 \\
    0 & AC
  \end{pmatrix},\, A\in \B_R,\, C\in\B_L\right\}.
\end{eqnarray*}
Soit $F=Diag(A,\,B)\in {G^0}^{\C}$, ($A,B\in SO(4)^{\C}$) alors on a
$\tau(F)=-L_{\flb[E]}FL_{\flb[E]}=Diag(B,\,A)$. Ainsi si $(\lambda \mapsto
F_{\lambda})\in \Lambda{G^0}^{\C}$, alors $F_{\lambda}\in\Lambda
{G_{\tau}^0}^{\C}$ \ssi en {\'e}crivant
$F_{\lambda}=\sum_{k\in\mathbb{Z}}\hat{F}_{2k}\lambda^{2k}$
on a $\hat{F}_{4k}\in\{Diag(R_aL_c,\,R_aL_c),\, a,c\in\h\otimes\C\}$
et $\hat{F}_{2k}\in\{Diag(R_aL_c,\,-R_aL_c),\, a,c\in\h\otimes\C\}$.
En particulier, si $(\lambda
\mapsto F_{\lambda})\in \Lambda_{\B}^{+}{G_{\tau}^0}^{\C}$ alors , comme
$F_0\in\B$, on en d{\'e}duit que $F_0\in\B_0$, donc
$\Lambda_{\B}^{+}{G_{\tau}^0}^{\C}=\Lambda_{\B_0}^{+}{G_{\tau}^0}^{\C}$.\\
On va utiliser le th{\'e}or{\`e}me suivant (cf. \cite{PS}) sur les groupes de
lacets.

\begin{theo}\label{iwasawa}
Soit $G$ un groupe de Lie compact, $G^{\C}$ son complexifi{\'e} et
$G^{\C}=G.\B_G$ sa d{\'e}composition d'Iwasawa. Alors
\begin{description}
  \item[(i)] La fonction produit
  $$
  \begin{array}{ccc}
    \Lambda G\times \Lambda_{\B_G}^{+}G^{\C} & \longrightarrow &
    \Lambda G^{\C} \\
    (\phi_{\lambda},\beta_{\lambda}) & \longmapsto &
    \phi_{\lambda}.\beta_{\lambda}
  \end{array}$$
est un diff{\'e}omorphisme.
  \item[(ii)] Il existe un ouvert $C_G$ de $\Lambda G^{\C}$ tel que
  la fonction produit
  $$
 \begin{array}{ccc}
\Lambda_{\ast}^{-}G^{\C}\times \Lambda^{+}G^{\C} & \longrightarrow & C_G \\
    (\phi_{\lambda}^{-},\phi_{\lambda}^{+}) & \longmapsto &
    \phi_{\lambda}^{-}.\phi_{\lambda}^{+}
  \end{array}$$
est un diff{\'e}omorphisme.
\end{description}
\end{theo}
\vspace{0.15cm}
\noindent On va en d{\'e}duire:

\begin{theo}\label{produit}
\begin{description}
  \item[(i)] La fonction produit
  $$
  \begin{array}{ccc}
    \Lambda G_{\tau}^0 \times \Lambda_{\B_0}^{+} {G_{\tau}^0}^{\C} &
    \longrightarrow &  \Lambda {G_{\tau}^0}^{\C} \\
    (F_{\lambda},B_{\lambda}) & \longmapsto &
    F_{\lambda}.B_{\lambda}
  \end{array}$$
est un diff{\'e}omorphisme.
  \item[(ii)] Il existe un ouvert $\overrightarrow{C}$ de
  $\Lambda {G_{\tau}^0}^{\C}$   tel que la fonction produit
  $$
 \begin{array}{ccc}
\Lambda_{\ast}^{-}{G_{\tau}^0}^{\C}\times \Lambda^{+}{G_{\tau}^0}^{\C} &
\longrightarrow & \overrightarrow{C} \\
    (F_{\lambda}^{-},F_{\lambda}^{+}) & \longmapsto &
    F_{\lambda}^{-}.F_{\lambda}^{+}
  \end{array}$$
est un diff{\'e}omorphisme.
\end{description}
\end{theo}

\noindent{\em D{\'e}monstration\hfil du\hfil th{\'e}or{\`e}me}\hfil \ref{produit}\hfil
 ---\hfil (i) \hfil
D'apr{\`e}s \hfil le \hfil th{\'e}or{\`e}me \hfil \ref{iwasawa}\,, \hfil l'application \\
$(F_{\lambda},B_{\lambda})\in  \Lambda G^0 \times \Lambda_{\B}^{+}
{G^0}^{\C} \longmapsto F_{\lambda}.B_{\lambda}\in \Lambda
{G^0}^{\C}$ est un diff{\'e}omorphisme. Soit $U_{\lambda}\in
\Lambda {G_{\tau}^0}^{\C}$ alors $U_{\lambda}=F_{\lambda}.B_{\lambda}$
avec $(F_{\lambda},B_{\lambda})\in \Lambda G^0 \times \Lambda_{\B}^{+}
{G^0}^{\C}$ et
 $F_{i\lambda}.B_{i\lambda}=\tau(F_{\lambda})\tau(B_{\lambda})$ or
$\tau(G^0)\subset G^0$ et $\tau(\B)\subset \B$ donc
$\tau(F_{\lambda})\in \Lambda G^0$,
$\tau(B_{\lambda})\in\Lambda_{\B}^{+} {G^0}^{\C}$ or d'apr{\`e}s le
th{\'e}or{\`e}me \ref{iwasawa}, il y a unicit{\'e} de la d{\'e}composition d'o{\`u}
$F_{i\lambda}=\tau(F_{\lambda}),B_{i\lambda}=\tau(B_{\lambda})$,
finalement $F_{\lambda}\in \Lambda G_{\tau}^0$, $B\in
 \Lambda_{\B}^{+} {G_{\tau}^0}^{\C}= \Lambda_{\B_0}^{+}{G_{\tau}^0}^{\C}$
 ce qui d{\'e}montre (i).\\
(ii) $\overrightarrow{C}=\Lambda_{\ast}^{-}{G_{\tau}^0}^{\C}.
\Lambda^{+}{G_{\tau}^0}^{\C}$ est l'image par un diff{\'e}omorphisme (celui donn{\'e}
par le th{\'e}or{\`e}me \ref{iwasawa}) de $\Lambda_{\ast}^{-}{G_{\tau}^0}^{\C}
\times \Lambda^{+}{G_{\tau}^0}^{\C}$ donc c'est une sous-vari{\'e}t{\'e} de
$\Lambda{G_{\tau}^0}^{\C}$ diff{\'e}omorphe {\`a} $\Lambda_{\ast}^{-}{G_{\tau}^0}^{\C}
\times \Lambda^{+}{G_{\tau}^0}^{\C}$. En outre $\overrightarrow{C}=
C_{G^0}\cap\Lambda {G_{\tau}^0}^{\C}$ donc c'est un ouvert de $\Lambda
{G_{\tau}^0}^{\C}$.\hfill $\blacksquare$ \\

Passons maintenant aux applications affines:

\begin{theo}\label{d{\'e}composition}
\begin{description}
  \item[(i)] On a la d{\'e}composition suivante:
  $$
\Lambda \G_{\tau}^{\C} = \Lambda \G_{\tau}. \Lambda_{\B_0}^{+}\G_{\tau}^{\C}
$$
  \item[(ii)]Il existe un ouvert $C$ de $\Lambda \G_{\tau}^{\C}$ tel
  qu'on ait la d{\'e}composition suivante:
  $$
 C =\Lambda_{\ast}^{-} \G_{\tau}^{\C}.\Lambda^{+}\G_{\tau}^{\C}
$$
\end{description}
\end{theo}

\noindent{\em D{\'e}monstration} ---
Il suffit de reprendre mot pour mot la d{\'e}monstration faite dans
\cite{HR1}, en remplacant $L_j$ par $L_{\flb[E]}$.\hfill $\blacksquare$ \\

\subsection{Repr{\'e}sentation de Weierstrass}
Nous allons suivre la m{\^e}me proc{\'e}dure que dans \cite{HR1} (i.e. utiliser
les m{\'e}thodes de \cite{DPW} en les adaptant).

\subsubsection{Potentiel holomorphe}

\begin{defi}
Soit $\mu$ une 1-forme sur $\Omega$ {\`a} valeurs dans $\Lambda
\gt_{\tau}^{\C}$, alors on dira de $\mu$ que c'est un potentiel
holomorphe si on a
$$
\mu_{\lambda}=\sum_{n\geq -2} \hat{\mu}_n \lambda^n
$$
avec $\hat{\mu}_n=\hat{\mu}_n(\frac{\partial}{\partial z}) dz$ o{\`u}
$\hat{\mu}_n(\frac{\partial}{\partial z})$ est holomorphe.
\end{defi}

\begin{theo}
Soit $U=(F,X)$ un rel{\`e}vement $\Sigma_V$ dont le $\rho_X$ est
harmonique et $U_{\lambda}=(F_{\lambda},X_{\lambda})\colon \Omega
\longrightarrow \Lambda\G_{\tau}$ son extension $\Sigma_V$ ($\Omega$ est
simplement connexe, on a choisi $z_0\in \Omega$ et $U_{\lambda}(z_0)=Id
\ \forall \lambda \in S^1$). Alors :
\begin{itemize}
  \item[$\bullet$]  Il existe une fonction holomorphe $H_{\lambda}\colon
  \Omega \to \Lambda \G_{\tau}^{\C}$ et une fonction
  $B_{\lambda}\colon \Omega \to
  \Lambda_{\B_0}^{+}\G_{\tau}^{\C}$ tel que $U_{\lambda}=
  H_{\lambda}.B_{\lambda}$.
  \item[$\bullet$]  En outre la forme de Maurer-Cartan $\mu_{\lambda}=
  H_{\lambda}^{-1}.dH_{\lambda}$ est un potentiel holomorphe: on dira que
  c'est un potentiel holomorphe pour $U_{\lambda}$ .
\end{itemize}
\end{theo}

\noindent{\em d{\'e}monstration} --- L'existence de $H_{\lambda}$ et
$B_{\lambda}$ est reli{\'e}e {\`a} la r{\'e}solution de l'{\'e}quation:
$$
0=\frac{\partial U_{\lambda}B_{\lambda}^{-1}}{\partial\bar{z}}=
U_{\lambda}\left(\alpha_{\lambda}\left(\frac{\partial}{\partial
\bar{z}}\right)-B_{\lambda}^{-1}\frac{\partial B_{\lambda}}{\partial
\bar{z}}\right)B_{\lambda}^{-1}
$$
qui est {\'e}quivalente {\`a}
$$
\frac{\partial B_{\lambda}}{\partial
\bar{z}}=B_{\lambda}\left(\alpha_0 + \lambda \alpha_1 +
\lambda^2\alpha_2\right)\left(\frac{\partial}{\partial\bar{z}}\right),
$$
avec la contrainte que $(\lambda\longmapsto
B_{\lambda}(z))\in\Lambda_{\B_0}^{+}\G_{\tau}^{\C}$ pour tout
$z\in\Omega$. L'existence de $B_{\lambda}$ s'obtient en suivant les
m{\^e}mes arguments que \cite{DPW}. Pour montrer que $\mu_{\lambda}$ est un
potentiel holomorphe, il suffit d'{\'e}crire:
$$
\mu_{\lambda}=H_{\lambda}^{-1}.dH_{\lambda}=B_{\lambda}(\alpha_{\lambda}-
B_{\lambda}^{-1}.dB_{\lambda}).B_{\lambda}^{-1}
$$
et d'utiliser le fait que $(\lambda\longmapsto
B_{\lambda}(z))\in\Lambda_{\B_0}^{+}\G_{\tau}^{\C}$ et que $z\mapsto
H_{\lambda}(z)$ est holomorphe.\hfill $\blacksquare$ \\

Inversement tout potentiel holomorphe produit une surface $\Sigma_V$
dont le $\rho_X$ est harmonique:

\begin{theo}\label{holomorphe}
Soit $\mu_{\lambda}$ un potentiel holomorphe, $z_0\in \Omega$ et
$H_{\lambda}^0$ une constante dans $\Lambda\G_{\tau}^{\C}$ (par
exemple $H_{\lambda}^0=Id$). Alors
\begin{itemize}
\item Il existe une unique fonction holomorphe $H_{\lambda}\colon
\Omega \longrightarrow \Lambda\G_{\tau}^{\C}$, tel que
$$
dH_{\lambda}=H_{\lambda}\mu_{\lambda} \quad et \quad
H_{\lambda}(z_0)=H_{\lambda}^0.
$$
\item En appliquant la d{\'e}composition
$\Lambda\G_{\tau}^{\C}=\Lambda\G_{\tau}.\Lambda_{\B_0}^{+}\G_{\tau}$
{\`a} $H_{\lambda}(z)$ pour chaque $z$, on obtient deux fonctions
$U_{\lambda}\colon \Omega\longrightarrow\Lambda\G_{\tau}$ et
$B_{\lambda}\colon\Omega\longrightarrow\Lambda_{\B_0}^{+}\G_{\tau}$ tel que
$H_{\lambda}(z)=U_{\lambda}(z).B_{\lambda}(z)\ \forall z\in\Omega$.
Alors $U_{\lambda}$ est une extension $\Sigma_V$ d'une surface
$\Sigma_V$ dont le $\rho_X$ est harmonique (pourvu que
$\hat{\alpha}_{-1}\neq 0$).
\item De plus $\mu_{\lambda}$ est un potentiel holomorphe pour
$U_{\lambda}$.
\end{itemize}
\end{theo}

\noindent{\em D{\'e}monstration} --- On a $d\mu_{\lambda}=0$, de plus
$\mu_{\lambda}(\frac{\partial}{\partial \bar{z}})=0$ implique
$\mu_{\lambda}\wedge\mu_{\lambda}=0$ d'o{\`u}
$$
d\mu_{\lambda} + \mu_{\lambda}\wedge\mu_{\lambda}=0,
$$
et il existe donc une unique fonction holomorphe
$H_{\lambda}\colon\Omega\longrightarrow\Lambda\G_{\tau}^{\C}$ tel que
\begin{equation}\label{Hmu}
dH_{\lambda}=H_{\lambda}.\mu_{\lambda}
\end{equation}
et $H_{\lambda}(z_0)=H_{\lambda}^0$. Ecrivons la d{\'e}composition du
th{\'e}or{\`e}me \ref{d{\'e}composition} pour $H_{\lambda}(z),\,z\in\Omega$:
il existe un unique couple
$(U_{\lambda},B_{\lambda})\in\Lambda\G_{\tau}\times\Lambda_{\B_0}^{+}
\G_{\tau}^{\C}$ tel que $H_{\lambda}(z)=U_{\lambda}(z).B_{\lambda}(z)$.
 Alors en utilisant (\ref{Hmu}), on a
\begin{equation}\label{U-1dU}
U_{\lambda}^{-1}.dU_{\lambda}=B_{\lambda}(\mu_{\lambda}-B_{\lambda}^{-1}.
dB_{\lambda})B_{\lambda}^{-1}
\end{equation}
Posons $\alpha=U_{\lambda}^{-1}.dU_{\lambda}$. Alors (\ref{U-1dU})
nous dit que $\alpha_{\lambda}$ doit s'{\'e}crire sous la forme
$$
\alpha_{\lambda}=\sum_{n\geq -2}\hat{\alpha}_n\lambda^n
$$
Mais comme $\alpha_{\lambda}$ est r{\'e}elle par d{\'e}finition, on a $\overline
{\hat{\alpha}}_n=\hat{\alpha}_{-n}$ et donc
$$
\alpha_{\lambda}=\lambda^{-2}\hat{\alpha}_{-2} + \lambda^{-1}\hat{\alpha}_{-1}
+\hat{\alpha}_0 + \lambda\hat{\alpha}_{1} + \lambda^{2}\hat{\alpha}_{2}.
$$
De plus, en utilisant $B_{\lambda}^{-1}=\hat{B}_0^{-1}-\lambda
\hat{B}_0^{-1}\hat{B}_1\hat{B}_0^{-1} + \cdots $, il r{\'e}sulte d'apr{\`e}s
(\ref{U-1dU}) que
$$
\hat{\alpha}_{-2}=\hat{B}_0\hat{\mu}_{-2}\hat{B}_0^{-1}\quad \text{et}
\quad \hat{\alpha}_{-1}=\hat{B}_0\hat{\mu}_{-1}\hat{B}_0^{-1} +
\hat{B}_1\hat{\mu}_{-2}\hat{B}_0^{-1} - \hat{B}_0\hat{\mu}_{-2}\hat{B}_0^{-1}
\hat{B}_1 \hat{B}_0^{-1},
$$
c'est \`{a} dire:
\begin{eqnarray*}
\hat{\alpha}_{-2} & = & \hat{B}_0\hat{\mu}_{-2}\hat{B}_0^{-1}\\
\hat{\alpha}_{-1} & = & \hat{B}_0\hat{\mu}_{-1}\hat{B}_0^{-1} +
[\hat{B}_1 \hat{B}_0^{-1},\hat{B}_0\hat{\mu}_{-2}\hat{B}_0^{-1}]\ .
\end{eqnarray*}
Ainsi on voit que $\hat{\alpha}_{-1}$ et $\hat{\alpha}_{-2}$ sont
des (1,0)-formes. De plus comme $\alpha_{\lambda}$ v{\'e}rifie
automatiquement la condition $d\alpha_{\lambda} +
\alpha_{\lambda}\wedge\alpha_{\lambda}=0$ alors le corollaire \ref{lacet}
implique -pourvu que $\hat{\alpha}_{-1}\neq 0$ - que $U_{\lambda}$ est une
extension $\Sigma_V$ dont le $\rho_X$ est harmonique.
Enfin comme $U_{\lambda}=H_{\lambda}.B_{\lambda}^{-1}$, on voit que
$\mu_{\lambda}$ est un potentiel holomorphe pour
$U_{\lambda}$.\hfill $\blacksquare$ \\

\subsection{Potentiel m{\'e}romorphe}

Le potentiel holomorphe construit au th{\'e}or{\`e}me \ref{holomorphe} est
loin d'{\^e}tre unique. On peut rem{\'e}dier {\`a} cela en incluant les
potentiels m{\'e}romorphes.

\begin{defi}
Un potentiel m{\'e}romorphe est une 1-forme m{\'e}romorphe sur $\Omega$ {\`a}
valeurs dans $\Lambda\gt_{\tau}^{\C}$ qui s'{\'e}crit
$$
\mu_{\lambda}=\lambda^{-2}\hat{\mu}_{-2} + \lambda^{-1}\hat{\mu}_{-1}.
$$
\end{defi}
\begin{theo}
Soit $U_{\lambda}:\Omega \longrightarrow\Lambda\G_{\tau}$ une
extension $\Sigma_V$ (dont le $\rho_X$ est harmonique). Alors il existe
une suite de points isol{\'e}s de $\Omega$ , $S=\{a_n,\,n\in I\}$ tel que :
\begin{itemize}
  \item $\forall z \in \Omega\smallsetminus\{a_n,\,n\in I\},\quad
  U_{\lambda}(z)\in
  C=\Lambda_{\ast}^{-}\G_{\tau}^{\C}.\Lambda^{+}\G_{\tau}^{\C}$.
  \item De plus les fonctions $U_{\lambda}^{-}\colon \Omega\smallsetminus
  \{a_n,\,n\in I\}\longrightarrow\Lambda_{\ast}^{-}\G_{\tau}^{\C} $ et
  $U_{\lambda}^{+}\colon \Omega\smallsetminus  \{a_n,\,n\in I\}
  \longrightarrow  \Lambda^{+}\G_{\tau}^{\C}$  ainsi d{\'e}finies sont
  holomorphes .
  \item En outre $U_{\lambda}^{-}$ s'{\'e}tend en un fonction m{\'e}romorphe
  sur $\Omega$.
  \item Sa forme de Maurer-Cartan $\mu_{\lambda}=(U_{\lambda}^{-})^{-1}
  dU_{\lambda}^{-}$ est un potentiel m{\'e}romorphe.
\end{itemize}
\end{theo}

\noindent{\em D{\'e}monstration} --- C'est une simple adaptation de la
d{\'e}monstration de \cite{DPW}.\hfill $\blacksquare$ \\

\begin{rem}{\em
Le potentiel holomorphe est unique par unicit{\'e} de la d{\'e}composition
$U_{\lambda}=U_{\lambda}^{-}.U_{\lambda}^{+}$. D'autre part, on peut
retrouver $X_{\lambda}$ en appliquant la m{\'e}thode du th{\'e}or{\`e}me
\ref{holomorphe} {\`a} $\mu_{\lambda}$. En effet, comme on peut
toujours supposer que $U_{\lambda}(z_0)=Id$, on a
$U_{\lambda}^{-}(z_0)=Id$ donc $U_{\lambda}^{-}$ est solution de
$\mu_{\lambda}=(U_{\lambda}^{-})^{-1}.dU_{\lambda}^{-},\,U_{\lambda}^{-}(z_0)
=Id.$ En outre, on peut {\'e}crire $U_{\lambda}^{+}=B_{\lambda}H$ avec
$H\in C^{\infty}(\Omega\smallsetminus S,G_0^0)$ et $B_{\lambda}\in
 C^{\infty}(\Omega\smallsetminus
 S,\Lambda_{\B_0}^{+}\G_{\tau}^{\C})$, alors on a
 $U_{\lambda}^{-}=U_{\lambda}(U_{\lambda}^{+})^{-1}=(U_{\lambda}H^{-1})
 B_{\lambda}^{-1}$ mais $U_{\lambda}H^{-1}$ est la partie $\Lambda \G_{\tau}$
dans la d{\'e}composition de $U_{\lambda}^{-}$ suivant
$\Lambda\G_{\tau}^{\C}=\Lambda\G_{\tau}.\Lambda_{\B_0}^{+}\G_{\tau}^{\C}$
sur $\Omega\smallsetminus S$ et $U_{\lambda}H^{-1}$ est un rel{\`e}vement de
$X_{\lambda}$.}
\end{rem}

\section{Le vecteur courbure moyenne}

On se propose de calculer le vecteur courbure moyenne d'une surface
$\Sigma_V$.\\
Soit $X\colon\Omega\subset\R^2 \to\R^8$ une immersion conforme $\Sigma_V$.
Alors on a
$$
H=\frac{e^{-2f}}{2}\,\triangle X=\frac{e^{-f}}{2}\left(\frac{\partial q}
{\partial u} + \frac{\partial q'}{\partial v} +f_u q +f_v q'\right).
$$
Posons $(E_1,E_2)=\mathcal{R}_{\rho}^{-1}(q,q')$, alors
$$
H=\frac{e^{-f}}{2}\mathcal{R}_{\rho}\left[ \frac{\partial
E_1}{\partial u}+ \frac{\partial E_2}{\partial v}
+  \begin{pmatrix}  0 & 0 \\    0 & R_{\gamma_u}  \end{pmatrix}E_1 +
\begin{pmatrix}  0 & 0 \\    0 & R_{\gamma_v}  \end{pmatrix}E_2 +
f_uE_1 +f_vE_2\right]
$$
avec comme d'habitude $\gamma=d\rho.\rho^{-1}$. En outre, en {\'e}crivant que
 $d(dX)=0$, on a
$$
\frac{\partial E_1}{\partial v} -\frac{\partial E_2}{\partial u} +
\begin{pmatrix}  0 & 0 \\    0 & R_{\gamma_v}  \end{pmatrix}E_1 -
\begin{pmatrix}  0 & 0 \\    0 & R_{\gamma_u}  \end{pmatrix}E_2 +
f_vE_1 -f_uE_2=0\ ,
$$
appliquons l'endomorphisme $L_{\flb[E]}$ {\`a} cette {\'e}quation
$$
\frac{\partial E_1}{\partial u}+ \frac{\partial E_2}{\partial v} +
  \begin{pmatrix}   -y.\gamma_v \\  0  \end{pmatrix} +
  \begin{pmatrix}   x.\gamma_u \\  0  \end{pmatrix} + f_vE_2 +f_uE_1
  =0.
$$
Injectons cela dans l'expression de $H$ :
\begin{eqnarray*}
H & = & \frac{e^{-f}}{2}\mathcal{R}_{\rho} \begin{pmatrix}
 y.\gamma_v -x.\gamma_u  \\ y.\gamma_u +x.\gamma_v  \end{pmatrix}\\
  & = & \frac{e^{-f}}{2}\mathcal{R}_{\rho}
\left[ \begin{pmatrix}
    -R_{\gamma_u} & 0 \\
    0 & R_{\gamma_u}
  \end{pmatrix}E_1   +
  \begin{pmatrix}
    -R_{\gamma_v} & 0 \\
    0 & R_{\gamma_v}
  \end{pmatrix}E_2 \right].
\end{eqnarray*}
En posant $\gamma^d=d\rho.\rho^{-1}$ et $\gamma^g=\rho^{-1}.d\rho$,
on peut r{\'e}{\'e}crire cela sous la forme
$$
H = \frac{e^{-f}}{2}\mathcal{R}_{\rho}\left[ \begin{pmatrix}
    -R_{\gamma_u^d} & 0 \\
    0 & R_{\gamma_u^g}
  \end{pmatrix}q   +
  \begin{pmatrix}
    -R_{\gamma_v^d} & 0 \\
    0 & R_{\gamma_v^g}
  \end{pmatrix}q' \right].
$$

\section{Surfaces $\omega_I$-isotropes, $\rho$-harmoniques}

\subsection{Le produit vectoriel de $\oct$ et le groupe $Spin(7)$}

\noindent Consid{\'e}rons l'application $ A\colon u\in \im\oct \mapsto
  \begin{pmatrix}
    L_u & 0 \\
    0 & -L_u
  \end{pmatrix}\in
\text{End}_{\R}(\oct)\oplus\text{End}_{\R}(\oct)$. On v{\'e}rifie que
$A(u)^2=-|u|^2Id$. Ainsi, $A$ est une application de Clifford et
on montre qu'elle se prolonge en un isomorphisme d'alg{\`e}bre $Cl(7)\cong
\text{End}_{\R}(\oct)\oplus\text{End}_{\R}(\oct)$ (cf. \cite{Har}).
On a de plus $Cl(7)^{\text{paire}}\cong\text{End}_{\R}(\oct)$ et
$Spin(7)$ est le groupe engendr{\'e} par $\{L_u,\, u\in S(\im\oct)\}$.
D'autre part, $g\in SO(8)$ appartient {\`a} $Spin(7)$ \ssi
$$
\forall u \in \im\oct=\R^7,\, \exists w\in \im\oct=\R^7 / $$
$$
g\,L_u\,g^{-1}=L_w
$$
\begin{sloppypar}
\noindent Cette propri{\'e}t{\'e} nous donne aussi la  repr{\'e}sentation
vectoriel de $Spin(7)$, \hbox{$\chi\colon Spin(7)\longrightarrow SO(7)$}
(rev{\^e}tement universel de $SO(7)$): $\forall g\in Spin(7),\,\forall u\in
\im\oct$
\begin{equation}\label{glug-1}
L_{\chi_g(u)}=g\,L_u\,g^{-1}
\end{equation}
ce qui donne $\chi_g(u)=g(ug^{-1}(1))$. On a alors que $g\in
O(\oct)$ est dans $Spin(7)$ \ssi
\begin{equation}\label{guv}
g(uv)=\chi_g(u)g(v)
\end{equation}
pour tout $u,v \in \oct$ (on pose $\chi_g(1)=1$).
\end{sloppypar}
Consid{\'e}rons maintenant le produit vectoriel de $\oct$:
$$
q\times q'=-\im(q.\bar{q'})=\im(q'.\bar{q})
$$
pour $q,q'\in \oct$. C'est une application antisym{\'e}trique  de
$\oct\times\oct$ dans $\im\oct$. Ainsi elle d{\'e}finit une application
$$
\begin{array}{crcl}
\rho \colon & Gr_2(\oct) & \longrightarrow & S(\im\oct)\\
            & q\wedge q' & \longmapsto  & q\times q'
\end{array}
$$
de la grassmanienne des plans orient{\'e}s de $\oct$ dans
$S^6 \subset \im\oct$. La propri{\'e}t{\'e} fondamentale qui va nous
permettre de faire dans un cadre plus g{\'e}n{\'e}ral ce que l'on a d{\'e}j{\`a}
fait pour les surfaces $\Sigma_V$ est la suivante
\begin{equation}\label{cross}
\forall g \in Spin(7),\ (g.q)\times (g.q')=\chi_g(q\times q').
\end{equation}
Elle veut dire que le produit vectoriel est $Spin(7)-${\'e}quivariant
(lorsque $Spin(7)$ agit sur $\oct\times\oct$ de mani{\`e}re naturelle et
sur $\im\oct=\R^7$ par l'interm{\'e}diaire de $\chi$). En particulier
lorsque l'on se restreint aux actions de $Spin(7)$ sur $Gr_2(\oct)$
et sur $S^6=S(\im\oct)$ alors $\rho$\, est $Spin(7)-${\'e}quivariante.\\
Pour tout couple orthonorm{\'e} de $\oct$, on a $\rho(q,q')=-q.\bar{q'}$
et donc $q'=\rho .q$. Ainsi l'ensemble des couples orthonorm{\'e}s de
$\oct$ est diff{\'e}omorphe {\`a} $S^6\times S^7$. De plus,  on rappelle que
l'on a pour tout $(q,q')\in \oct^2$
$$
q\times q'=-\sum_{i=1}^7 \omega_i(q,q')e_i \ .
$$
Consid{\'e}rons  maintenant, pour $I\varsubsetneqq\{1,2,\ldots 7\}$, les ensembles
\begin{eqnarray*}
V_I & = & \{(q,q')\in S^7\times S^7/ \langle q,q'\rangle =\omega_i(q,q')=0,\, i\in I\}\\
Q_I & = & \{P\in Gr_2(\oct) /\ \omega_i(P)=0,\, i\in I\}.
\end{eqnarray*}
On a $Q_I=V_I/SO(2)$; et $\rho(Q_I)=S^{I}=S(\bigoplus_{i\notin I}\R
e_i)$. En particulier  $V_I\cong S^7\times S^{6-|I|}$. Pour
$I=\varnothing$, on a $Q_{\varnothing}=Gr_2(\oct)$. Pour $I=\{1,2,3\}$
on retrouve l'ensemble $Q$ {\'e}tudi{\'e} dans la section~2.\par
Cherchons le sous-groupe, $G_I$, de $Spin(7)$ qui conserve  les
$\omega_i$, $i\in I$. $g\in Spin(7)$ conserve $\omega_i$ \ssi il
commute avec $L_{e_i}$, $gL_{e_i}g^{-1}=L_{e_i}$, ce qui veut dire
que $\chi_g(e_i)=e_i$, c'est {\`a} dire qu'il conserve $e_i$ par son
action sur $S^6$. Il en r{\'e}sulte que
$G_I=\chi^{-1}(SO(\bigoplus_{i\notin I}\R
e_i))\simeq\chi^{-1}(SO(7-|I|))=Spin(7-|I|)$. Ainsi $G_I$ agit  sur
$Q_I$ ainsi que que sur $S^{I}$ et le produit vectoriel $ \rho\colon Q_I
\to S^{I}$ est $G_I-${\'e}quivariant. De plus, $G_I$ agit transitivement
sur $S^{I}$. En outre le stabilisateur d'un point de $S^6$ pour
l'action de $Spin(7)$ s'identifie {\`a} $\chi^{-1}(SO(6))=Spin(6)$ et
donc $S^6=Spin(7)/Spin(6)$. Plus g{\'e}n{\'e}ralement on a
$S^{I}=G_I/G_{I\cup\{k\}}\simeq Spin(7-|I|)/Spin(6-|I|)$.

On peut toujours consid{\'e}rer en toute g{\'e}n{\'e}ralit{\'e} que $I=\{1,2...,|I|\}$. En
effet, le groupe $SO(7)$ agit transitivement sur les $S^{I}$ avec
$I$ de m{\^e}me cardinal. Donc en prenant l'image r{\'e}ciproque  par
$\chi$, et $\rho$, on voit que $Spin(7)$ agit transitivement sur les
$V_I$ avec $I$ de cardinal fix{\'e}. Ainsi quitte {\`a} faire agir un
{\'e}l{\'e}ment de $Spin(7)$ on peut toujours consid{\'e}rer que
$I=\{1,2,...,|I|\}$, $S^{I}=S^{6-|I|}$, et $G_{I}=Spin(7-|I|)$.

\begin{exem} Prenons $I=\{1,2,3\}$, on retrouve le cas {\'e}tudi{\'e} dans
les sections pr{\'e}c{\'e}dentes. On a $Q_I=Q$ , $G_I=Spin(4)=S^3\times
S^3$ et $S^3\simeq S^3\times S^3/S^3=Spin(4)/Spin(3)$. Plus
pr{\'e}cis{\'e}ment, on a vu que $G_I=\left\{
  \begin{pmatrix}
   R_a & 0 \\
    0 & R_b
  \end{pmatrix},\,a,b\in S^3\right\}$ et $Stab_{G_I}(\flb[E])=G_{I\cup\{4\}}
=\left\{\begin{pmatrix}
   R_a & 0 \\
    0 & R_a
  \end{pmatrix},\,a\in S^3\right\}$. En outre, on a vu que l'action de
  $G_I$ sur $S^3$ s'{\'e}crit
$$
\chi_g(u)=\bar{a}u\,b  \qquad i.e.
$$
$$
\chi\colon Diag(R_a,R_b)\in Spin(4)\mapsto L_{\bar{a}}R_b\in SO(4)
$$
\end{exem}

\noindent Nous allons maintenant passer en revue les autres cas possibles.
\begin{enumerate}
  \item [0.] Si $I=\varnothing$, alors $G_I=Spin(7)$ et
  $V_I=St_2(\oct)$ l'ensemble des couples orthonorm{\'e}s de $\oct$.
  L'action de $Spin(7)$ sur $S^6$ est donn{\'e}e par
  $\chi_g(u)=g(u)\overline{g(1)}$ d'apr{\`e}s (\ref{guv}). $Spin(7)$
  agit transitivement sur les couples orthonorm{\'e}s de $\oct$. En
  effet, l'action de $Spin(7)$ sur $S^7$ est transitive et
  $Stab(1)=G_2$, et donc $S^7=Spin(7)/G_2$. Ensuite l'action de
  $G_2$ sur $S(\{1\}^{\perp})=S^6$ est transitive et
  $Stab_{G_2}(e_1)=SU(\bigoplus_{i>1}\R e_i,\,L_{e_1})\simeq SU(3)$
  et donc $S^6=G_2/SU(3)$. On d{\'e}duit alors facilement de cela que $Spin(7)$
   agit transitivement  sur $St_2(\oct)$ et $Stab_{Spin(7)}(1,e_1)=SU(3)$ et
   donc  $St_2(\oct)=Spin(7)/SU(3)$, $Gr_2(\oct)=Spin(7)/(SO(2)\times
   SU(3))$.
  \item  Si $I=\{1\}$ alors $G_I=Spin(6)=SU(4)$. En effet, c'est le
  sous-groupe de $Spin(7)$ qui commute avec la structure complexe
  $L_{e_1}$, il est donc inclus dans $U(4)$. D'autre part comme
  $Spin(6)$ est engendr{\'e} , dans $Cl(6)$, par les $u\cdot v$, $(u,v)$ couple
  orthonorm{\'e}  de $\R^6$, alors comme $(u\cdot v)^{2}=-1$, si on note $A$
   la repr{\'e}sentation de $Spin(6)$ dans $\text{End}(\C^4)$ alors on a
   $A(u\cdot v)^2=-Id$, donc $A(u\cdot v)$ est une structure
   complexe de $\C^4$ donc $\det_{\C}(A(u\cdot v))\in\{\pm 1,\pm
   i\}$, ainsi $\det_{\C}(A(Spin(6)))\subset\{\pm 1,\pm i\}$. Enfin la
   connexit{\'e} de $Spin(6)$ et un calcul des dimensions nous donne
   $A(Spin(6))=SU(4)$.\\
   En outre, $Spin(6)=SU(4)$ agit transitivement sur $V_I$ qui n'est
   autre que l'ensemble des couples hermitiens de $(\oct,L_{e_1})=\C^4$, et
   donc il agit transitivement sur $Q_I$. Le stabilisateur d'un
   couple $(q,q')$ est isomorphe {\`a} $SU(2)$. Ainsi
   $V_{\{1\}}=SU(4)/SU(2)$, $Q_{\{1\}}=SU(4)/(SO(2)\times SU(2))$.
   \item Si $I=\{1,2\}$, alors $G_I=Spin(5)\simeq U(2,\h)$. En
   effet, c'est le sous-groupe de $Spin(7)$ qui commute avec
   $I_1=L_{\flh[I]}$ et $I_2=L_{\flh[J]}$ donc avec $I_3=I_1I_2$
   aussi. Il est donc inclus dans le groupe des automorphismes $\h
   -$lin{\'e}aires pour la structure de $\h -$espace vectoriel sur
   $\oct$ d{\'e}finie par $I_1,I_2,I_3$; mais il est aussi inclus dans
   $SO(8)$, il est donc inclus dans le groupe des $\h
   -$automorphismes de $\oct$ qui conservent la forme hermitienne
   quaternionique $C=\langle \cdot ,\cdot\rangle + \langle \cdot,I_1
   \cdot\rangle i + \langle\cdot,I_2 \cdot\rangle j + \langle\cdot ,I_3
   \cdot\rangle k$,
   qui est un conjugu{\'e} de $U(2,\h)$. Ensuite un calcul des
   dimensions et la connexit{\'e} des deux groupes permettent de conclure
   que $Spin(5)\simeq U(2,\h)$.
   $Spin(5)$ n'agit pas  transitivement sur $V_{\{1,2\}}$ car
   $Spin(5)$ conserve  la forme $C$ et agit librement et
   transitivement sur les fibres de $C$, tandis que $C(V_{\{1,2\}})=
   [-1,1]k$. Ainsi les orbites sont caract{\'e}ris{\'e}es par la fonction
   $C$ qui sur $V_{\{1,2\}}$ vaut $C(q,q')=\langle q,I_3.q'\rangle k $
    et chaque orbite est diff{\'e}omorphe {\`a} $Spin(5)$.
   \item Le cas $I=\{1,2,3\}$ a d{\'e}j{\`a} {\'e}t{\'e} longuement {\'e}tudi{\'e} dans les
   sections pr{\'e}c{\'e}dentes. On l'a rappel{\'e} dans l'exemple~1. En outre
   rajoutons que d'apr{\`e}s la section~2 (cf. la d{\'e}monstration du
   th{\'e}or{\`e}me~\ref{p}), le groupe $Spin(4)$ n'agit pas transitivement sur $Q$,
    et que les orbites sont caract{\'e}ris{\'e}es par la fonction
   $r\colon q\wedge q'\mapsto Im(x.\bar{x'})$. Il y a alors une orbite de
   dimension 6\,: $G.P_1=r^{-1}(0)$, une $S^2-$famille d'orbites de dimension
   5, dont la r{\'e}union est $G.P_2$: $\{P/r(P)=u/2\}$, $u$ d{\'e}crivant $S^2$,
   et enfin une famille d'orbite de dimension 6, celles tel que
   $0<|r|<1/2$.
   \item Si $|I|=4$, on a vu que $G_I=Spin(3)=\{Diag(R_a,R_a),\,a\in
   S^3\}$. L'action de $Spin(3)$ sur $S^2$ s'{\'e}crit
   $\chi_g(\rho)=\bar{a}\rho\,a$ i.e $\chi\colon Diag(R_a,R_a)\in
   Spin(3)\mapsto L_{\bar{a}}R_a\in SO(3)$.
   \item Si $|I|=5$ alors
   $G_I=Spin(2)=\{Diag(R_{e^{i\theta}},R_{e^{i\theta}}),\,\theta\in
   \R\}$. L'action de $Spin(2)$ sur $S^1$ est donn{\'e}e par
   $$
   \chi_g(e^{i\beta})=e^{-2i\theta}.e^{i\beta}
   $$
   l'identification entre $\C$ et $\R e_6\oplus\R e_7$ {\'e}tant $ae_6
   +be_7\mapsto a+ib$.
   \item Si $|I|=6$, alors $G_I=\{\pm Id\}$, $S_I=\{\pm e_7\}$ et
   $V_I=\{(q,\pm L_{e_7}q),\, q\in S^7\}$.
\end{enumerate}
\vspace{1.5mm}
On a besoin, pour refaire dans le cas g{\'e}n{\'e}ral ce que l'on a fait
pr{\'e}c{\'e}demment, de d{\'e}finir un application $\tilde{\rho}\colon G_I\to
S^{I}$. Soit donc $I\varsubsetneqq\{1,...7\}$ et $e\in S(\im\oct)
\smallsetminus\{e_i,\,i\in I\}$, disons pour que les choses soient
fix{\'e}es une fois pour toutes que $e=e_{|I|+1}$ (en supposant comme on
en a le droit que $I=\{1,...,|I|\}$). Posons alors
$$
\tilde{\rho}_I(g)=\chi_g(e)
$$
pour tout $g\in G_I$, alors $\tilde{\rho}_I(G_I)=S^{I}$ et par passage au
quotient $\tilde{\rho}_I$ d{\'e}finit l'isomorphisme
$G_{I}/G_{I\cup\{k\}}\simeq S^{I}$. De plus on a
$$
\rho(g^{-1}.q,g^{-1}.q')=e\Longleftrightarrow \tilde{\rho}_I(g)=\rho(q,q').
$$
En outre $\rho(q,q')=e\Longleftrightarrow q'=L_e q$.

\subsection{Alg{\`e}bres de Lie}

L'automorphisme int{\'e}rieur de $Spin(7)$, $int_{L_e}$, stabilise
$G_I$. En effet, soit $g\in G_I$, alors pour $i\in I$ on a
$$
\begin{array}{cclcl}
  L_{e_i}(L_e\,g\,L_e^{-1})L_{e_i}^{-1} & = & L_{e_i}(L_e\,g\,L_e)L_{e_i}
   & = & (-L_eL_{e_i})\,g\,(-L_{e_i}L_e) \\
   & = & L_e(L_{e_i}\,g\,L_{e_i}^{-1})L_e^{-1} & = & L_e\,g\,L_e^{-1}
\end{array}
$$
ainsi $L_e\,g\,L_e^{-1}$ commute avec $L_{e_i},\,i\in I$ d'o{\`u}
$L_e(G_I)L_e^{-1}\subset G_I$. L'alg{\`e}bre de Lie de $G_I$,
$\g_I=spin(7-|I|)$ est stable par $Ad\,L_e$ et on a $(Ad\,L_e)^2=Id$ donc
elle se d{\'e}compose en somme directe des deux sous-espaces propres
associ{\'e}s aux valeurs propres $\pm 1$ respectivement:
$$
\g_I=\g_0(I)\oplus\g_2(I)
$$
avec
$$
\g_0(I)=\ker(Ad\,L_e -Id)\,, \quad \g_2(I)=\ker(Ad\,L_e + Id).
$$
On posera $\g_k=\g_k(\varnothing)$, $k=0,2$, ainsi $spin(7)=\g_0\oplus
\g_2$. Pour $|I|=5$, on a
$int_{L_e}=-Id$ d'o{\`u} $\g_2(I)=0$. Pour $|I|=6$, on a $\g_I=0$.\\
Consid{\'e}rons maintenant le groupe $\G=\mathcal{A}Spin(7)$ des
isom{\'e}tries affines de $\oct$ dont la partie lin{\'e}aire est dans
$Spin(7)$ et plus g{\'e}n{\'e}ralement $\G_I=\mathcal{A}G_I$ et $\g(I)$ son
alg{\`e}bre de Lie. Soit $\tau_e$ l'automorphisme int{\'e}rieur de $\G_I$
d{\'e}fini par $(-L_e,0)$. Il induit un automorphisme d'ordre 4 de
$\g(I)$ (il est d'ordre 4 pour tout $|I|$ m{\^e}me pour $|I|=5$ et 6,
car sur $\oct$, $L_e$ est d'ordre 4) ce qui donne une d{\'e}composition
de $\g^{\C}(I)=\g(I)\otimes\C$ suivant les espaces propres de $\tau_e$
associ{\'e}s aux valeurs propres $i^k,\,-1\leq k\leq 2$. On notera
$\g_k^{\C}(I)$ ces espaces propres. Les espaces propres associ{\'e}s aux
valeurs propres $\pm 1$ ne d{\'e}pendent pas de $I$, on les notera
simplement $\g_{\pm 1}^{\C}$. On a alors
$$
  \begin{array}{lclclcl}
   \g_{-1}^{\C} & = & \ker(L_e - iId) & , & \g_1^{\C} & = & \ker(L_e + iId) \\
    \g_0^{\C}(I) & = & \g_0(I)\otimes\C & , & \g_2^{\C}(I) & = & \g_2(I)
    \otimes\C\
  \end{array}
$$
On a $[\g_k^{\C}(I),\,\g_l^{\C}(I)]\subset
\g_{k+l\text{ mod\,4}}^{\C}(I)$, puisque $\tau_e$ est un
automorphisme.

\subsection{Surfaces $\omega_I-$isotropes et $\rho -$harmoniques}

\begin{defi}
Soit $\Sigma$ une surface immerg{\'e}e de $\oct$, alors il lui est
associ{\'e}e une application $\rho_{\Sigma}\colon \Sigma \to S^6$
d{\'e}finie par $\rho_{\Sigma}(z)=\rho(T_z\Sigma)$ i.e. si
$X\colon\Sigma\to\oct$ est une immersion alors
$\rho_{\Sigma}=X^{\ast}\rho$. On dira que $\Sigma$ est $\rho
-$harmonique si $\rho_{\Sigma}$ est harmonique. Soit
$I\varsubsetneqq \{1,...,7\}$ alors on dira que $\Sigma$ est
$\omega_I-$isotrope si $\forall z \in  \Sigma,\ T_z\Sigma\in Q_I$
(si $I=\varnothing$ il n'y a aucune condition).
Dans ce cas, $\rho_{\Sigma}$ est {\`a} valeurs dans $S^{I}\subset S^6$.
\end{defi}
\begin{defi} On appellera rel{\`e}vement $\omega _I-$isotrope (si
$I=\varnothing$ on dira seulement rel{\`e}vement) une
application $U=(F,X)\colon \Sigma\to \G_I$ tel que $X$ soit une
immersion conforme $\omega_I-$isotrope et $\tilde{\rho}_I\circ F
=\rho_{\Sigma}$
\end{defi}
\begin{theo}
Soit $\Omega$ un ouvert simplement connexe, $\alpha \in
T^{\ast}\Omega\otimes\g(I)$, alors
\begin{itemize}
  \item [$\bullet$] $\alpha$ est la forme de Maurer-Cartan d'un
  rel{\`e}vement $\omega_I-$isotrope \ssi
  $$
d\alpha +\alpha\wedge\alpha=0,\quad \alpha_{-1}''=0\ \text{et }\,
\alpha_{-1}' \text{ ne s'annule pas,}
$$
  \item[$\bullet$] dans ce cas, $\alpha$ correspond {\`a} une immersion
  conforme $\omega_I-$isotrope, $\rho -harmonique$ \ssi la forme de
  Maurer-Cartan prolong{\'e}e $\alpha_{\lambda}=\lambda^{-2}\alpha_2' +
  \lambda^{-1}\alpha_{-1}' + \alpha_0 +
  \lambda\alpha_1''+\lambda^2\alpha_2''$ v{\'e}rifie
  $$
  d\alpha_{\lambda} +\alpha_{\lambda}\wedge\alpha_{\lambda}=0,\
  \forall\lambda\in \C^{\ast}.
  $$
  \end{itemize}
\end{theo}

\noindent \emph{D{\'e}monstration} --- Pour le premier point cf. la
section~3. Pour le second on remarque que l'on a
$$
 d\alpha_{\lambda} +\alpha_{\lambda}\wedge\alpha_{\lambda}=
  d\beta_{\lambda^2} +\beta_{\lambda^2}\wedge\beta_{\lambda^2}
$$
o{\`u} $\beta_{\lambda}=\lambda^{-1}\alpha_{2}' + \alpha_0 + \lambda\alpha_2''$
est la forme de Maurer-Cartan prolong{\'e}e associ{\'e}e {\`a} $\beta=F^{-1}.dF$,
la forme de Maurer-Cartan du rel{\`e}vement $F\in G_I$ de $\rho_X\in
S^{I}$. D'apr{\`e}s \cite{DPW} (cf. aussi \cite{H1}, ou \cite{H2}) on
sait que $\rho_X$ est harmonique \ssi $d\beta_{\lambda} +\beta_{\lambda}
\wedge\beta_{\lambda}=0\ \forall \lambda \in S^1$ ce qui ach{\`e}ve la
d{\'e}monstration.\hfill $\blacksquare$ \\

On peut maintenant faire ce que l'on a fait dans la section~4
{\`a} l'aide des groupes de lacets et obtenir une repr{\'e}sentation de
Weierstrass par des potentiels holomorphes {\`a} valeurs dans
$\Lambda\g^{\C}(I)_{\tau}$.

\begin{rem}{\em
On a choisi $e=e_{|I|+1}$. Avec ce choix le groupe de sym{\'e}trie
$G_{I\cup\{|I|+1\}}$ n'est autre que le groupe d'isotropie pour
l'action de $G_I$. Donc c'est pratique lorsqu'on fait une {\'e}tude
th{\'e}orique o{\`u} on envisage successivement les diff{\'e}rents cas
possibles. Mais ce choix n'est pas pratique si l'on veut que les
d{\'e}compositions des alg{\`e}bres de Lie se corespondent i.e. la
d{\'e}composition du sous-cas soit la "trace"  de la d{\'e}composition du cas
plus g{\'e}n{\'e}ral. Pour cela il faut donc prendre le m{\^e}me $e$. Il suffit
de le prendre compatible avec le sous-cas  alors il sera compatible
avec le cas plus g{\'e}n{\'e}ral. Exemple: Si on veut {\'e}tudier le cas $|I|=2$
ainsi que le sous-cas $|I|=3$, alors le choix de $e=\flb[E]$ pour
les deux convient bien et les d{\'e}compositions se correspondent:
$\g_k^{\C}(\{1,2,3\})=\g_k^{\C}(\{1,2\})\cap \g^{\C}(\{1,2,3\})$. Donc on peut
faire une {\'e}tude simultan{\'e}e des deux (on a le m{\^e}me automorphisme
$\tau_e$) et la repr{\'e}sentation de Weierstrass du sous-cas
s'obtiendra tout simplement en prenant dans la repr{\'e}sentation de Weierstrass
du cas plus g{\'e}n{\'e}ral, des potentiels holomorphes {\`a} valeurs dans une sous
alg{\`e}bre de Lie, $\Lambda\g^{\C}(J)_{\tau}$, de l'alg{\`e}bre de Lie du cas plus
 g{\'e}n{\'e}ral, $\Lambda\g^{\C}(I)_{\tau}$.}
\end{rem}
\begin{rem}{\em
On peut {\'e}videmment envisager le cas plus g{\'e}n{\'e}ral des surfaces
$\omega_E-$isotropes o{\`u} $E$ est un sous-espace vectoriel de
$\im\oct$ : on dira qu'un plan $P$ de $\oct$ est $\omega_E-$isotrope
s'il est isotrope pour $\omega_e=\langle\cdot,L_e\cdot\rangle$ pour tout
$e\in S(E)$.
Pour tout $E$, il existe $g\in Spin(7)$ qui envoie l'ensemble des
plans $\omega_E-$isotropes sur l'ensemble des plans
$\omega_I-$isotropes avec $|I|=\dim E$. En effet, soit $h\in SO(7)$
tel que $h.E=\text{Vect}(e_i,\,i\in I)$ alors pour $g\in\chi^{-1}(\{h\})$ on a
$g^{\ast}\omega_I=\omega_E$ avec
$\omega_E=\text{Vect}(\omega_e,\,e\in E)$, $\omega_I=\text{Vect}
(\omega_i,\,i\in I)$. Donc quitte {\`a} faire agir un {\'e}l{\'e}ment fixe de $Spin(7)$
on est ramen{\'e} au cas que l'on a {\'e}tudi{\'e}.}
\end{rem}
\begin{rem}{\em\begin{sloppypar}
Ce que l'on vient de faire dans $\oct$ peut {\^e}tre fait dans $\h$. On
d{\'e}finit le produit vectoriel $x\times y=-\im(x.\bar{y})$ pour
$x,y\in\h$. On a $x\times y =-\sum_{i=1}^{3}\omega_i(x,y)e_i$ avec
$\omega_i=\langle\cdot,e_{i}\,\cdot\rangle$ et $(e_1,e_2,e_3)=(i,j,k)$ la base canonique de
$\im\h=\R^3$. Alors on a pour $g=R_aL_b\in SO(4)$, $(gx)\times
(gy)=-\im(bx\bar{y}\bar{b})=b(x\times y)b^{-1}$. Ainsi la
repr{\'e}sentation vectorielle de $Spin(3)=\{L_b,\,b\in S^3\}$,
\hbox{$\chi\colon Spin(3)\to SO(3)$} est donn{\'e}e par $L_b\in Spin(3)\mapsto
int_b=L_bR_{b^{-1}}\in SO(\im\h)$. En proc{\'e}dant comme on l'a fait
dans $\oct$, on montre que les surfaces de $\R^4$ $\rho
-$harmoniques sont un syst{\`e}me compl{\`e}tement int{\'e}grable. Plus
g{\'e}n{\'e}ralement, les surfaces $\omega_I-\text{isotropes}$, $\rho -$harmoniques
de $\R^4$ sont un syst{\`e}me compl{\`e}tement int{\'e}grable. Ici on a
$|I|=0,1$ ou 2. Pour $|I|=1$ on retrouve les surfaces lagrangiennes
hamiltoniennes stationnaires. Pour $|I|=2$ on trouve les surfaces sp{\'e}ciales
lagrangiennes.\\
D'ailleurs une surface $\omega_I-$isotrope, $\rho -$harmonique de
$\h$ n'est autre qu'une surface $\omega_I-$isotrope, $\rho -$harmonique
de $\oct$ contenue dans le sous-espace $\h$ de $\oct$. D'autre part,
on voit que si l'immersion $X$ est {\`a} valeurs dans $\R^3=\im\h$ alors
$X$ est $\rho -$harmonique \ssi elle est {\`a} courbure moyenne
constante. Ainsi l'ensemble des CMC de $\R^3$ n'est autre que
l'ensemble des surfaces $\rho -$harmoniques de $\h$ incluses dans
$\im\h$.
\end{sloppypar}}
\end{rem}

\subsection{Calcul du vecteur courbure moyenne}
Dans le cas des surfaces $\Sigma_V$ ($|I|=3$) on disposait du
rel{\`e}vement $\mathcal{R}\colon\rho\in S^3\mapsto
\mathcal{R}_{\rho}\in Spin(4)$, en particulier la fibration $Spin(4)\to S^3$
est triviale (on a vu que c'est un produit semi-direct). Ceci nous
a permis d'{\'e}crire l'{\'e}quation lin{\'e}aire (\ref{dedzbar}) qui caract{\'e}rise
les surfaces $\Sigma_V$ et de calculer le vecteur courbure moyenne
en utilisant le fait que le fibr{\'e} $Q\to S^3$ est trivial: on a
repr{\'e}sent{\'e} $(q,q')$ par $(\rho,(E_1,E_2))$. Dans le cas g{\'e}n{\'e}ral, ce
n'est pas possible: $SO(n+1)\to SO(n+1)/SO(n)=S^n$ n'est pas trivialisable
sinon la sph{\`e}re $S^n$ serait parall{\'e}lisable or on sait que ce n'est le cas
que pour $n=1,3,7$. On retrouve donc le cas $n=3$ et le cas {\'e}vident $n=1$
(dans ce cas $Spin(2)=S^1$). Dans le cas g{\'e}n{\'e}ral donc on ne peut pas
faire la s{\'e}paration pr{\'e}c{\'e}dente. Cependant, on peut toujours calculer
le vecteur courbure moyenne en fonction de $\rho$ et la formule
obtenue est valable pour n'importe quelle surface de $\oct$ sans
aucune hypoth{\`e}se.\\
Soit $X\colon\Omega\to \oct$ une immersion conforme alors par d{\'e}finition
de $\rho\colon Gr_2(\oct)\to \im\oct$ on a
$$
\ast dX=-\rho_X.dX
$$
et cette {\'e}quation d{\'e}termine le $\rho_X=X^{\ast}\rho$ de l'immersion
(i.e. si on a $\ast dX=-\sigma.dX$ alors $\sigma=\rho_X$). Prenons
la diff{\'e}rentielle de cette {\'e}quation, on obtient
$$
\triangle X=(\partial_v \rho)\cdot{\frac{\partial X}{\partial u}}
-(\partial_u\rho)\cdot{\frac{\partial X}{\partial v}}=\rho\cdot\left(
(\partial_u\rho)\cdot{\frac{\partial X}{\partial u}} + (\partial_v \rho)
\cdot{\frac{\partial X}{\partial v}}\right)
$$
On en d{\'e}duit
\begin{eqnarray*}
H=\frac{e^{-2f}}{2}\triangle X  & = & \frac{e^{-2f}}{2}\left((\partial_v \rho)
\cdot{\frac{\partial X}{\partial u}}
-(\partial_u\rho)\cdot{\frac{\partial X}{\partial v}}\right)\\
  & = & \frac{e^{-2f}}{2}\,\rho\cdot\left((\partial_u\rho)\cdot
{\frac{\partial X}{\partial u}} +(\partial_v \rho)
\cdot{\frac{\partial X}{\partial v}}\right).
\end{eqnarray*}

\vspace{0.2cm}
\noindent \textbf{Idrisse Khemar\\
khemar@math.jussieu.fr}

\end{document}